\newif\ifsmfart
\IfFileExists{smfart.cls}
   {\documentclass[12pt,english]{smfart}
    \IfFileExists{smfenum.sty}{\usepackage{smfenum}}{}
    \usepackage{bull}
    \smfarttrue}
{\message{^^J*** It would be better to typeset this file with smfart.cls ***^^J^
^J}

\documentclass[12pt]{amsart}}
\usepackage{amssymb}
\usepackage{epsfig}
\usepackage{times}
\numberwithin{equation}{section}

\input xy
\xyoption{all}

\makeatother
\makeatletter

\def\Title    {Reconstruction of function fields}
\def\Author   {Fedor Bogomolov and Yuri Tschinkel}
\def\Subject  {Algebraic geometry}
\def\Keywords {Galois groups, function fields}

\newif\ifpdf
   \ifx\pdfoutput\undefined
   \pdffalse
   \else
           \pdfoutput=1
            \pdftrue
   \fi
      \ifpdf
\usepackage[pdftex,
hyperindex=true,
pdftitle={\Title},pdfauthor={\Author},pdfsubject={\Subject},
pdfkeywords={\Keywords},pdfpagemode={UseOutlines},
bookmarks=true,bookmarksopen=true,
bookmarksopenlevel=0,bookmarksnumbered=true,
pdfstartview={FitV},
colorlinks=true
]{hyperref}
 \else
 \usepackage{hyperref}
 \fi

\usepackage{graphicx}

\advance\headheight 2pt
\calclayout
\allowdisplaybreaks[3]

\theoremstyle{plain}
\newtheorem{prop}[subsection]{Proposition}
\newtheorem{Theo}[section]{Theorem}
\newtheorem{thm}[subsection]{Theorem}
\newtheorem{coro}[subsection]{Corollary}

\newtheorem{lemm}[subsection]{Lemma}

\newtheorem{defn}[subsection]{Definition}

\theoremstyle{definition}

\theoremstyle{remark}
\newtheorem{rem}[subsection]{Remark}

\newtheorem{exam}[subsection]{Example}

\newtheorem{nota}[subsection]{Notations}

\def\G{{\mathcal G}}   % abel. pro-l-Galois group 
\def\D{{\mathcal D}}   % decomposition group
\def\I{{\mathcal I}}   % inertia groups

\def\Res{{\varrho}}
\def\pic{{\varphi}}
\def\Ga{{\Gamma}}%\value group
\def\Val{{\mathcal V}}
\def\DVal{{\mathcal D}{\mathcal V}}

\def\KK{\boldsymbol{K}}

\def\sK{{\mathfrak K}}

\def\AF{{\mathcal F}}
\def\LF{{\mathcal L}}
\def\NS{{\rm NS}}

\def\FS{{\mathcal F}{\mathcal S}}

\def\supp{{\rm supp}}
\def\Ext{{\rm Ext}}

\def\no{\noindent}
\def\rk{{\rm rk}}
\newcommand{\trdeg}{{\rm tr}\, {\rm deg}}
\newcommand{\ovl}{\overline}

\newcommand{\La}{\Lambda}

\newcommand{\al}{\alpha}

\def\cH{{\mathcal H}}

\def\cL{{\mathcal L}}
\def\cM{{\mathcal M}}

\def\dv{{\rm div}}
\newcommand{\Hom}{{\rm Hom}}
\newcommand{\Ker}{{\rm Ker}}

\def\lra{\longrightarrow}
\def\ra{\rightarrow}

\def\F{{\mathbb F}}

\def\P{{\mathbb P}}
\def\Q{{\mathbb Q}}

\def\Z{{\mathbb Z}}

\def\N{{\mathbb N}}

\def\sc{{\mathsf c}}

\def\bA{{\mathbb A}}

\def\bP{{\mathbb P}}

\def\mc{{\mathfrak c}}

\def\mo{{\mathfrak o}}

\def\mm{{\mathfrak m}}
\def\ml{{\mathfrak l}}
\def\mL{{\mathfrak L}}

\def\mP{{\mathfrak P}}

\def\Pic{{\rm Pic}}
\def\Div{{\rm Div}}

\def\pr{{\rm pr}}

\def\char{{\rm char}}

\author{Fedor Bogomolov}
\address{Courant Institute of Mathematical Sciences, N.Y.U. \\
 251 Mercer str. \\
 New York, NY 10012, U.S.A.}
\email{bogomolo@cims.nyu.edu}

\author{Yuri Tschinkel}
\address{Department of Mathematics \\
         Princeton University\\
         Fine Hall, Washington Road\\
         Princeton, NJ 08544-1000,  U.S.A.}
\email{ytschink@math.princeton.edu}        

\title[Function fields]{Reconstruction of function fields}

\begin{document}

\begin{abstract}
We study the structure of abelian subgroups of
Galois groups of function fields of surfaces.
\end{abstract}
\date{\today}

\maketitle
\tableofcontents

\setcounter{section}{0}       
\section*{Introduction}
\label{sect:introduction}

We fix two primes $p$ and $\ell$. We will assume that 
$\ell\neq p$ and $p\neq 2$.
Let $k=\ovl{\F}_p$ be an algebraic closure of
the finite field $\F_p$. 
Let $X$ be an algebraic variety defined
over $k$ and $K=k(X)$ its function field.
We will refer to $X$ as a {\em model} of $K$.  
Let $\Pic(X)$ be the Picard group, $\NS(X)$
the N\'eron-Severi group of $X$ and 
$\G^a_K$ the abelianization of 
the pro-$\ell$-quotient $\G_K$ of the
absolute Galois group of $K$. 
Under our assumptions on $k$, $\G^a_K$ 
is a torsion free $\Z_{\ell}$-module.  
Let  $\G^c_K$ be its canonical
central extension - the second lower 
central series quotient of $\G_K$.
It determines the following structure 
on $\G_K^a$: a set $\Sigma_K$ of
distinguished (primitive) subgroups which are isomorphic
to {\em finite rank} 
(torsion free) $\Z_{\ell}$-modules.
A topologically noncyclic subgroup $\sigma\in \Sigma_K$ iff
\begin{itemize}
\item $\sigma$ lifts to an abelian subgroup of $\G^c_K$;
\item $\sigma$ is maximal:
there are no abelian subgroups $\sigma'\subset \G^a_K$
which lift to an abelian subgroup of $\G^c_K$
and contain $\sigma$ as a proper subgroup.
\end{itemize}
We will call $\Sigma_K$ a fan.
The main theorem of this paper is

\begin{Theo}                                     
\label{thm:main}
Let $K$ and $L$ be function fields over
algebraic closures of finite fields
of characteristic $\neq 2,\ell$. 
Assume that $K=k(X)$ is a function field of a surface $X/k$ 
such that
\begin{itemize}
\item[(1)] $\Pic(X)=\NS(X)$;
\item[(2)] there exists an isomorphism
$$
\Phi=\Phi_{K,L}\,:\, \G^a_K\simeq \G^a_{L}
$$
of abelian pro-$\ell$-groups inducing a bijection of sets
$$
\Sigma_K = \Sigma_{L}.
$$
\end{itemize}
Then $L$ is a finite purely inseparable extension of $K$.
\end{Theo}

\no
The conditions on $X$ do not depend on a choice of a model, which 
we can assume to be smooth (by resolution of singularities 
in dimension two). 
We chose to treat in detail a class of surfaces for which the proof 
of Theorem~\ref{thm:main} is most transparent. 
The assumption (1) is not necessary; we have included it since 
it allows us to avoid certain geometric technicalities.

In this paper we implement the program outlined in 
\cite{B-1} and \cite{B-3} describing the correspondence between 
higher-dimensional function fields and their abelianized Galois
groups.
For results concerning the reconstruction of 
function fields from their (full) Galois groups
(the birational Grothendieck program)
we refer to the works of Pop, Mochizuki and Efrat
(see \cite{pop}, \cite{mochizuki} and \cite{efrat}). \

\

\no
{\bf Acknowledgments.} 
Both authors were partially supported by the NSF. 
The second author was employed by the Clay Mathematics Institute. 
We are grateful to Laurent Lafforgue and Barry Mazur for their interest. 
We thank the participants of the 
Algebraic Geometry Seminar at the University of Nice
for their comments and suggestions.

\section{Basic algebra and geometry of fields}
\label{sect:basicalg}

In this section we collect some technical results about 
function fields.

\begin{nota}
\label{nota:k}
Throughout, $k$ is an algebraic closure of the finite field $\F=\F_p$
and $K=k(X)$ the function field of an algebraic variety
$X/k$ over $k$ (its {\em model}).
\end{nota}

We assume familiarity with
\begin{itemize}
\item basic notions of field theory (transcendence degree, 
purely inseparable extensions);
\item basic notions of algebraic geometry
(models $X$ of a field $K$, 
$k$-rational points $X(k)$,
Picard group $\Pic(X)$, 
N\'eron-Severi group $\NS(X)$). 
\end{itemize}

\begin{lemm}
\label{lemm:C}
Let $C/k$ be a smooth curve and $Q\subset C(k)$ a finite set.
Then there exists an $n=n_Q\in \N$ such that for every
degree zero divisor $D$ with support in $Q$ the divisor 
$nD$ is principal.  
\end{lemm}

\begin{proof}
Every finitely generated subgroup of a torsion group
is finite. Since the group of degree zero divisors $\Pic^0(C)$ 
(over any finite field) is a torsion group and its subgroup
of divisors with support in a finite set of points
$Q\subset C(k)$ is finitely generated, the claim follows.
\end{proof}

\begin{lemm}
\label{lemm:cc}
Let $K/k$ be the function field of a surface,
$C/k$ a smooth curve on a model of $K$ and 
$Q=\{ q_0,...,q_s\}\subset C(k)$ 
a finite set of points. Then there exist a model 
$X$ of $K$, irreducible divisors 
$D_j$, $H_j,H_j'$ on $X$ (with $j=0,...,s$)
and a positive integer $n=n_Q$ 
such that:
\begin{itemize}
\item[(1)] $X$ is smooth and contains $C$;
\item[(2)] $D_j\cap C=q_j$ for all $j=1,...,s$;
\item[(3)] $n(D_j-D_0)$ restricted to $C$ is a principal divisor;
\item[(4)] $n(D_j-D_0) +(H_j-H_j')$ is a principal divisor on $X$; 
\item[(5)] the divisors $D_j$ are pairwise disjoint;
\item[(6)] all intersections between $D_j, H_i$ and 
$H_i'$ are transversal, pairwise
distinct and outside $C$;
\item[(7)] $H_j,H_j'$ don't intersect $C$. 
\end{itemize}
\end{lemm}

\begin{proof}
On a model $X$ containing $C$ as a smooth curve
choose any divisors $D_j\subset X$ 
passing (transversally) through $q_j$ (for all $j=0,...,s$). 
Blowing up points in $C(k)\setminus Q$ 
we can insure that the (strict transform of) $C$ becomes contractible
and that the image of the surface
under a contracting morphism is {\em projective}. 
Blowing up again (if necessary) and removing components
of exceptional divisors, we can insure that the (strict transforms) 
$D_j\cap C=q_j$ (for all $j$). 
By Lemma~\ref{lemm:C}, there exists an $n=n_Q$ such that the restriction of
$n(D_j-D_0)$ to $C$ is a principal divisor. 
We continue to blow up (outside $Q$)
so that each $n(D_j-D_0)$ becomes  a trivial line bundle on some affine
neighborhood of $C$ in some model $X$. 
Throughout, $C$ remains contractible
and we write 
$$
\pi\,:\, X\ra Y
$$
for the corresponding blow-down.  
Now $n(D_j-D_0)$ is induced from a line bundle on 
$Y$ (which is projective). In particular,
there exist {\em ample} classes  $[H_j],[H_j']\in \Pic(Y)$ 
such that
$$
[n(D_j-D_0)] +([H_j]-[H_j'])
$$
is a principal divisor on $X$ (here we identified $[H_j],[H_j']$
with their full transforms in $X$). 
Finally, we can choose representatives $H_j,H_j'\subset Y$
of these classes which are disjoint from $\pi(C)$, 
irreducible and satisfy all required transversality assumptions. 
\end{proof}

\begin{lemm}
\label{lemm:purely}
Let $K/\sK$ be a purely inseparable extension. Then 
\begin{itemize}
\item $\sK\supset k$;
\item $K/\sK $ is a finite extension;
\item $\sK=k(X')$ for some algebraic variety $X'$.
\end{itemize}
\end{lemm}

\begin{defn}
\label{defn:gener}
We write $\ovl{E}^{K}\subset K$ for the {\em normal} 
closure of a subfield $E\subset K$ 
(elements in $K$ which are algebraic over $E$). 
We say that $x\in K\setminus k$ 
is {\em generating} if $\ovl{k(x)}^{K}=k(x)$. 
\end{defn}

\begin{rem}
\label{rem:just}
If  $E\subset K$ is 1-dimensional 
then for all $x\in E\setminus k$ one has
$\ovl{k(x)}^{K}=\ovl{E}^{K}$ (a finite extension of $E$).  
\end{rem}

\begin{lemm}
\label{lemm:subf}
For any subfield $E\subset K$ there is a canonical
sequence  
$$
X\stackrel{\pi_E}{\lra} Y'\stackrel{\rho_E}{\lra} Y,
$$
where 
\begin{itemize}
\item $\pi_E$ is birational dominant with irreducible generic fiber;
\item $\rho_E$ is quasi-finite and dominant; 
\item $k(Y')=\ovl{E}^{K}$ and $k(Y)=E$.
\end{itemize}
\end{lemm}

For generating $x\in K$ we write 
$$
\pi_x\,:\, X\ra Y
$$
for the morphism from 
Lemma~\ref{lemm:subf}, with $k(Y)=k(x)$.  
For $y\in K\setminus k(x)$  we define  
$\deg_x(y)$ (the degree of $y$ on the generic
fiber of $\pi_x$) as the degree of the corresponding 
surjective map from the generic fiber 
of $\pi_x$ under the projection $\pi_y$.

\begin{prop}
\label{prop:geome}
Let $X/k$ be a smooth surface. Then 
\begin{itemize}
\item[(1)] if $X$ contains finitely many rational curves then
the same holds for every 
model $X'$ of $K$;
\item[(2)] if $\Pic(X)=\NS(X)$  
then for every 1-dimensional subfield $E\subset K=k(X)$
such that $E=\ovl{E}^{K}$ one has $E= k(x)$ for some $x\in K$;
\item[(3)] for every curve $C\subset X$ and every finite set of 
irreducible divisors
$D_1,...,D_s$ of $X$ not containing $C$, there exists a blowup
$$
\pi\,:\, \tilde{X}\ra X
$$
such that every branch in the strict transform $\pi^{-1}(C)$ 
intersects at most one irreducible component of the full transform 
of $\cup_{j=1}^s D_j$ and these intersections are transversal. 
\end{itemize}
\end{prop}

\begin{proof}
Property (1) is evident. 
Property (2) follows from the fact that
every such 1-dimensional field 
corresponds to a dominant map $X\ra C$ onto a curve. 
If $\Pic(X)=\NS(X)$ then $\Pic^0(X)=0$ and $X$ admits 
no such maps onto curves of genus $\ge 1$.  
The last property follows from resolution of singularities 
for surfaces.
\end{proof}

\begin{lemm}
\label{lemm:XX}
Assume that $\Pic(X)=\NS(X)$ and let 
$x,y\in k(X)\setminus k$ be such that 
$$
\deg_x(y)=\min_{f\in K\setminus \ovl{k(x)}^{K}}(\deg_x(f)).
$$
Then $y$ is generating: $k(y)= \ovl{k(y)}^{K}$.
\end{lemm}

\begin{proof}
If $y$ is not generating then $y=z(y')$ for some $y'\in K=k(X)$
and some function $z\in k(y')$ of degree $\ge 2$. 
This implies that $\deg_x(y)\ge 2\deg_x(y')$, 
contradicting minimality. 
\end{proof}

\begin{prop}
\label{prop:gener}
Let $X/k$ be an algebraic variety of dimension $\ge 2$ 
such that $\Pic(X)=\NS(X)$. 
If $t\in K=k(X)$ is not generating 
then there exist $y,y'\in K$ such that
for all $\kappa_1,\kappa_2,\kappa_1',\kappa_2'\in k^*$ 
the elements
\begin{equation}
\label{eqn:fib}
y, y/(t+\kappa_1), 
(y+\kappa_2)/t \text{ and } y', y'/(t+\kappa_1'), (y'+\kappa_2')
\end{equation}
are generating and the elements
\begin{equation}
\label{eqn:fibb}
1,y,y',t
\end{equation}
are linearly independent over $k$.
\end{prop}

\begin{proof}
Write $E:=\ovl{k(t)}^{K}$.  
By Proposition~\ref{prop:geome} (2), there exists an  
$x\in K$ such that $E=k(x)$ (so that $t=t(x)$). 
We have a dominant morphism
$$
\pi_x\,:\, X\ra \P^1
$$
with irreducible generic fiber.
Choose two algebraically independent $y$ and $y'$ so that 
$$
\deg_x(y)=\deg_x(y')=\min_{f\in K\setminus \ovl{k(x)}^{K}}(\deg_x(f))
$$
and 
$$
1,y,y',t
$$
are linearly independent 
(linear independence can be checked, for example, 
by restricting to a fiber of $\pi_x$). 
By Lemma~\ref{lemm:XX}, both $y$ and $y'$ are generating. It suffices to
observe that all elements in \eqref{eqn:fib}
have the same degree on the generic fiber of $\pi_x$.  
\end{proof}

The next proposition is a characterization of multiplicative
groups of subfields $\sK\subset K$. We will say that $y\in K^*$ 
is a {\em power} if there exist an 
$x\in K^*$ and an integer $n\ge 2$ such that
$y=x^n$. 

\begin{prop}
\label{prop:KKK}
Let $X$ be an algebraic variety of dimension $\ge 2$
such that every 1-dimensional subfield $E\subset K=k(X)$ 
has the form $E=k(x)$ for some $x\in K^*$. 
Let $\sK^*\subset K^*$ be a subset such that
\begin{itemize}
\item[(1)] $\sK^*$  is closed under multiplication;
\item[(2)] for every $E=k(x)\subset K$ with $E=\ovl{E}^K$ 
there exists an $r=r(x)\in \N$ 
such that 
$$
\sK^*\cap E^*=(E^*)^{r};
$$
\item[(3)] there exists a 
$y\in K\setminus k$ with $r(y)=1$ such that $y$ is
not a power.   
\end{itemize}
Then $\sK:=\sK^*\cup 0$ is a field and $K/\sK$ is a purely inseparable
finite extension. 
\end{prop}

\begin{proof}
Once we know that $\sK$ is a field we can 
conclude that every 
$x\in K^*$ is either in $\sK^*$ or some power of it
is in $\sK^*$. Of course, it can only be 
a power of $p$ so that $K/\sK$ is 
a purely inseparable extension, of finite degree (by  
Lemma~\ref{lemm:purely}).

By (3), $k\subset \sK$. To conclude that $\sK$ is a field, 
it suffices to show
that for every $x\in \sK$ 
one has $x+1\in \sK$ (and then use multiplicativity). 
For every $x\in \sK\setminus k$ 
with $r(x)=1$ we have $\sK^*\cap k(x)^*=k(x)^*$ and 
$$
x+\kappa \in \sK^*, \, \text{ for all } \kappa \in k.
$$ 
In particular, this holds for $y$. 

Now consider $x\in \sK^*$ with $r(x)>1$. We claim that there exists a $\kappa\in k$
such that 
$$
z:=\frac{x+y+\kappa}{y+\kappa-1} \in \sK \,\, \text{ and } r(z)=1.
$$
This implies that
$$
z-1=(x+1)/(y+\kappa-1)\in \sK^*\,\, \text{  and } x+1\in \sK^*,
$$ 
(by multiplicativity).

To prove the claim,  choose a model $X$ of $K$ and 
consider the morphisms
$$
\begin{array}{ccc}
\pi_x\,: X &\ra & \P^1=(x:1) \\
\pi_y\,: X & \ra & \P^1=(y:1)
\end{array}
$$ 
(as in Lemma~\ref{lemm:subf}).
Since $x$ and $y$ are algebraically independent ($r(x)>1$), 
only finitely many components of the fibers of $\pi_x$ are contained in 
the fibers of $\pi_y$ and there exists a 
$\kappa\in k$ such that both fibers 
$$
\pi_y^{-1}(-\kappa) \text{ and } \pi_y^{-1}(1-\kappa)
$$
are transversal to the fibers of $\pi_x$. 

Then 
$$
\dv_0(y+\kappa-1)\not\subset \dv(x+y+\kappa),
$$
since $y+\kappa=-1$ on $y+\kappa-1$ and $x$ is nonconstant on these fibers
(where $\dv_0$ is the divisor of zeroes). 
It follows that {\em both} 
$$
t:=(y+\kappa)/x \,\, \text{ and } \,\, z:=(x+y+\kappa)/(y+\kappa-1)
$$
are not powers so that $r(t)=r(z)=1$.  
To show that $z\in \sK$ observe
that both  $x,y+\kappa\in \sK$ so that $t\in \sK$. 
Therefore, 
$$
t+1=(x+y+\kappa)/x\in \sK
$$ 
and, by (1), $x+y+\kappa\in \sK$. 
Finally, since  $(y+\kappa-1)\in \sK$ we get $z\in \sK$. 
\end{proof}

\section{Projective structures}
\label{sect:proj-str}

In this section we explain the connection between 
fields and axiomatic projective geometry. We follow closely
the exposition in \cite{mihalek72}.

\begin{defn}
\label{defi:proj}
Let $S$ be a (nonempty) set and $\mL=\mL(S)$ a collection of subsets 
${\mathfrak l}\subset S$ such that
\begin{itemize}
\item[P1] there exist an $s\in S$ and an $\mathfrak l\in \mL$
such that $s\notin \mathfrak l$;
\item[P2] for every $\mathfrak l\in \mL$ there exist at least
three distinct $s,s',s''\in \mathfrak l$;
\item[P3] for every pair of distinct $s,s'\in S$ there exists exactly 
one 
$$
\mathfrak l={\mathfrak l}(s,s')\in \mL
$$ 
such that $s,s'\in \mathfrak l$;
\item[P4] for every quadruple of pairwise distinct
$s,s',t,t'\in S$ one has
$$
{\mathfrak l}(s,s')\cap {\mathfrak l}(t,t')\neq \emptyset\,\, \Rightarrow\,\, 
{\mathfrak l}(s,t)\cap {\mathfrak l}(s',t')\neq \emptyset  . 
$$ 
\end{itemize}
Such a pair $(S,\mL)$ is called a {\em projective structure} on 
$S$ and the elements $\mathfrak l\in \mL$ are called {\em lines}.
\end{defn}

For $s\in S$ and $S'\subset S$ define the {\em join}
$$
s\vee S':=
\{ s''\in S\,|\, s''\in \ml(s,s')\, \text{ for some } s'\in
S'\} .
$$
For any finite set of points $s_1, \ldots , s_n$ define
$$
\langle s_1,\ldots s_n\rangle := s_1\vee \langle s_2 \vee \cdots \vee s_n\rangle
$$ 
(this does not depend on the order of
the points). Write $\langle S'\rangle $ for the join
of a {\em finite} set $S'\subset S$.
A finite set $ S'\subset S$ of pairwise distinct points 
is called {\em independent} 
if for all $s'\in S'$ one has
$$ 
s'\notin \langle S'\setminus \{ s'\}\rangle .
$$
A set of points $S'\subset S$ {\em spans} a set of points $T\subset S$ if 
\begin{itemize}
\item $\langle S''\rangle \subset T$ for every finite set $S''\subset S'$;
\item for every $t\in T$ there exists 
a finite set of points $S_t\subset S'$
such that $t\in \langle S_t\rangle$. 
\end{itemize} 
A set $T\subset S$ spanned by an 
independent set $S'$ of points of cardinality $\ge 2$ is called a
projective {\em subspace} of dimension $|S'|-1$.

A projective structure $(S,\mL)$ satisfies {\em Pappus' axiom}
if for all 2-dimensional subspaces and every configuration
of six points and lines in these subspaces as below

\

\centerline{\includegraphics[width=.5\textwidth]{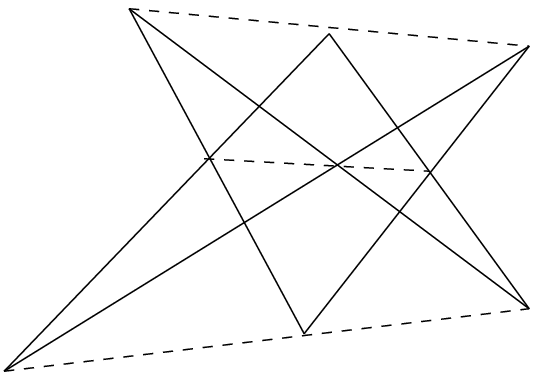}}

\

\noindent
the intersections are collinear. 
The main theorem of abstract projective geometry is:

\begin{thm}
\label{thm:abstr}
Let $(S,\mL)$ be a projective structure of dimension $n\ge 2$ 
which satisfies Pappus' axiom. 
Then there exists a field $k$ such that 
$S=\P_k^n$. This field is unique, up to isomorphism.
\end{thm}

\begin{proof}
See \cite{mihalek72}, Chapter 6. 
\end{proof}

\begin{defn}
\label{defi:inject}
A morphism of projective structures 
$$
\rho\,:\, (S,\mL)\ra (S',\mL')
$$ 
is an injection of sets $\rho\,:\, S\hookrightarrow S'$ 
such that $\rho(\ml)\in \mL'$ for all $\ml\in \mL$.
\end{defn}

\begin{exam}
\label{exam:Kk}
Let $k$ be a field and $\P^n_k$ the usual projective 
space over $k$ of dimension $n\ge 2$. 
Then $\P^n_k$ carries a natural projective structure: the set of lines
is the set of usual projective lines $\P^1_k\subset \P^n_k$.
Every (not necessarily finite) extension of fields 
$K/k$ induces a morphism of projective structures:
$\P^n_k\ra \P^n_K$.
\end{exam}

\begin{defn}
\label{defn:compa}
Let $S$ be an abelian group.
A projective structure $(S,\mL)$ on $S$ 
will be called {\em compatible} with 
the group structure if for all $s\in S$ and $\mathfrak l\in \mL$
one has  
$$
s\cdot {\mathfrak l} \in {\mathfrak L}(S).
$$ 
\end{defn}

\begin{exam}
\label{exam:basic}
Let $K/k$ be a field extension (not necessarily finite). Then 
$$
S:=\P_k(K)=(K\setminus 0)/k^*
$$ 
carries a natural projective structure
which is compatible with multiplication in the abelian group
$K^*/k^*$.
\end{exam}

\begin{thm}
\label{thm:skk}
Let $K/k$ and $K'/k'$ be field extensions of degree $\ge 4$ and 
$$
\bar{\phi}\,:\, S=\P_k(K)\ra \P_{k'}(K')=S'
$$ 
a bijection of sets 
which is an isomorphism of abelian groups 
and of projective structures.
Then 
$$
k\simeq k'\,\, \text{ and } \,\, K\simeq K'.
$$  
\end{thm}

\begin{proof}
Choose a plane $\P^2\subset S$ containing the 
identity $e\in S$, and two lines $\ml_1,\ml_2$ in this plane 
passing through $e$.
The set of all points $\P^2\setminus \{ \ml_1,\ml_2\}$
is a principal homogeneous space under the group
of projective automorphisms of $\P^1_k $ ($=\ml_1$) 
stabilizing one point (the intersection $\ml_1\cap \ml_2$). 
A choice of an additional point $s\in \mathfrak P$ outside 
$\ml_1\cup \ml_2$ trivializes this homogeneous space to the group 
of affine transformations of an affine line over $k$. 
In particular, it determines both
the additive and the multiplicative structure on $k$. 
This implies that $k$ is isomorphic to $k'$ and that for every
finite-dimensional space $V\subset K$ there exists
a unique $k'$-linear space $V'\subset K'$
such that the map $\bar{\phi}_V\,:\, \P_k(V)\ra \P_{k'}(V')$ 
lifts to a $(k,k')$-linear map $\phi_V\,:\, V\ra V'$.
Such a lift is unique modulo multiplication by a nonzero scalar in $k$ 
on the left (resp. $k'$ on the right). 
We can identify $V$ with $\P(V)\times k^*\cup \{0\}$ (as a set).  
If $V$ is such that $e\in \P(V)$ then there is a unique lift  
$\phi_V$ with the property $\bar{\phi}_V(e)=e'\in S'$.

Let $x,y\in K\setminus k$ be any elements
projecting to $\bar{x},\bar{y}\in \P_k(K)$ and  
$V\subset K$ a $k$-vector subspace containing
$$
1,x,y,xy. 
$$ 
Fix $\phi=\phi_V$ as above. 
Since $\bar{\phi}$ is an isomorphism of abelian groups there is a
$c(x,y)\in k^* $ such that
$$
\phi(x\cdot y)= 
\phi(x)\phi(y)c(x,y).
$$ 
We need to show that $c(x,y)=1$. For any $a\in k^*$ we have
$$ 
\phi((a+ x)\cdot y)
= \phi(a\cdot y + c({x},{y})\cdot x\cdot y) =
a'\cdot y' + c'({x},{y})\cdot x'\cdot y'   \in V'\subset K',
$$
by $(k,k')$-linearity of $\phi$. Since $\bar{\phi}$ preserves products, 
the right side must be $k'$-proportional to 
$$
a'\cdot y'+x'\cdot y'.
$$
On the other hand, $y'$ and $x'\cdot y'$ are $k'$-linearly independent
(since $x'\notin k'$).  This implies that $c'(x',y')=1$, as claimed. 
\end{proof}

\begin{defn}
\label{defi:partt}
A {\em partial projective structure} is a pair
$(S,\mP)$, where $\mP$ is a set of subsets of $S$ (lines)
such that for every triple of pairwise distinct 
elements $r,s,t\in S$ there exist distinct elements
$$
x,y,x',y'\in S
$$
(pairwise distinct from $r,s,t$)
and lines
$$
\ml(y,r),\ml(y,s),\ml(t,x),\ml(y',r),\ml(y',s),\ml(t,x'),
\ml(y,y') 
\in \mP
$$
with the property that
$$
\begin{array}{l}
r,x,y\in \ml(y,r), r,x',y'\in \ml(y',r),    \\
y,s \in \ml(y,s),y',s \in \ml(y',s), \\
t,x\in \ml(t,x), t,x'\in \ml(t,x'),\\
\ml(t,x)\cap \ml(y,s)\neq \emptyset, \ml(t,x')\cap \ml(y',s)\neq \emptyset,\\
\ml(y,y')\cap \ml(t,x)=\ml(y,y')\cap  \ml(t,x')=\emptyset.
\end{array}
$$
\end{defn}

\

\centerline{\includegraphics[width=.5\textwidth]{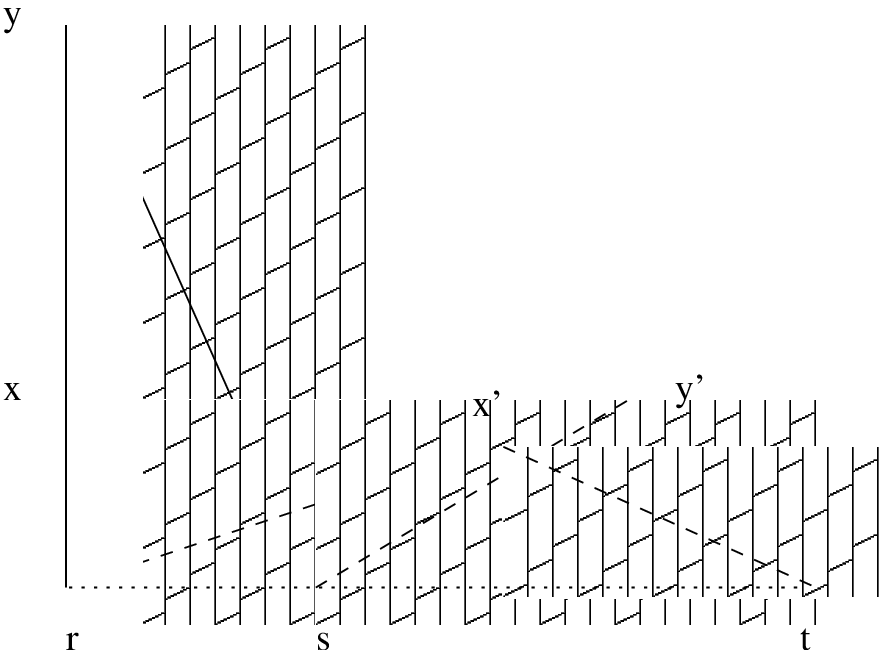}}

\

\begin{rem}
\label{rem:ppr}
Of course, a projective structure on a projective space of dimension $\ge 3$
is also a partial projective structure.
\end{rem}

\begin{exam}
\label{exam:basic2}
Let $K/k$ be the function field of an algebraic variety $X$ 
of dimension $\ge 2$. We have seen in Example~\ref{exam:basic}
that $S:=\P_k(K)$ carries a projective structure compatible with 
multiplication. Let $\mP$ be the set of all lines 
passing through $1$ and a generating element of 
$K$ (see Definition~\ref{defn:gener}) and all their translates by
elements in $S$ (under multiplication). 
Proposition~\ref{prop:gener} 
implies that $(S,\mP)$ is a partial projective structure.
Indeed, (after translation) we can assume that $r=1,s=t+\kappa$ etc.
The lines in $\mP$ containing 1 will be called {\em  primary}. 
\end{exam}

\begin{prop}
\label{prop:unique-p}
Let $(S,\mL)$ and $(S,\mL')$ be two 
projective structures on $S$ of dimension $\ge 3$ 
and assume that the intersection 
$\mL\cap \mL'$ contains a subset $\mP$ such that $(S,\mP)$ is 
a partial projective structure on $S$. Then $\mL=\mL'$.  
\end{prop}

\begin{proof}
For any tuple of points $(r,s)$ we need to show that 
the line $\ml(r,s)\in \mL$ is also in $\mL'$ (and vice versa). 
Let $t$ be any point distinct from $r$ and $s$. 
There exist points $x,y,x',y'$ as in Definition~\ref{defi:partt}
together with corresponding lines in $\mP$. Moreover,  
$t\in \ml(r,s)$ iff $t$ lies both in the plane spanned
by $\ml(y,r)$ and $\ml(y,s)$ and in the plane spanned by
$\ml(y',r)$ and $\ml(y',s)$ (which are distinct by the
assumption that $\ml(y,y')$ intersects neither $\ml(t,x)$  
not $\ml(t,x')$). These conditions are satisfied by 
the assumptions on the lines $\ml(t,x)$, resp. $\ml(t,x')$.
Thus $t\in \ml(r,s) $ iff $t\in\ml'(r,s)$. 
\end{proof}

\section{Flag maps}
\label{sect:af}

\begin{nota}
\label{nota:pr}
We fix two prime numbers $\ell$ and $p\neq 2,\ell$.
Let
\begin{itemize}
\item  $\F=\F_p$ be the finite field with $p$ elements;
\item $\F^*=\F\setminus \{ 0\}$ its multiplicative group;
\item $A$ a vector space over $\F$ of dimension $\dim(A)\in \N\cup \infty$;
\item $\bP(A)=(A\setminus 0)/\F^*$ its projectivization;
\item $\cM(A,S)$ the set of maps from $A$ to a set $S$;
\item for $\mu\in \cM(A,S)$ and $B\subset A$, 
$\mu_B$ the restriction of $\mu$ to $B$.
\end{itemize}
\end{nota}

\begin{defn}
\label{defn:aa}
A (complete) {\em flag} on a finite-dimensional $A$ is 
a collection of subspaces $(A_n)_{n=0,..,\dim(A)}$ 
such that
\begin{itemize}
\item $A_0=A$;
\item $A_n\supsetneq A_{n+1}$ for all $n=0,...,\dim(A)-1$.
\end{itemize}
In particular, 
$A_{n}\setminus A_{n+1} \neq \emptyset$, for all $n$,
and $A_{\dim(A)} =\{ 0\}$.
\end{defn}

\begin{defn}
\label{defn:homo}
A map $\mu\in \cM(A,S)$ will be called homogeneous
if for all $a\in A$ and all $\kappa\in \F^*$ one 
has 
$$
\mu(\kappa\cdot a)=\mu(a).
$$
\end{defn}

We think of homogeneous maps as being defined on 
the projectivization $\bP(A)$ and write $\cM(\bP(A),S)$ for
the space of such maps.

\begin{defn}
\label{defn:inv}
A map $\mu\in \cM(A,S)$ will be 
called a flag map if
\begin{itemize}
\item $\mu$ is homogeneous;
\item every  finite-dimensional $\F$-vector space $B\subset A$
has a complete flag $(B_n)_{n=0,...,\dim(B)}$  such that
$\mu$ is constant on $B_n\setminus B_{n+1}$, for all
$n=0,...,\dim(B)-1$.  
\end{itemize}
The set of flag maps will be denoted by
$\AF(A,S)$ or $\AF(\bP(A),S)$.
\end{defn}

\begin{rem}
The flag property does not 
depend on the value of $\mu$ in $0\in A$. 
Thus we will generally work on $A\setminus 0$ and $\bP(A)$. 
\end{rem}

\begin{defn}
\label{defn:lf}
Let $S$ be a ring and $A$ an $\F$-algebra. 
A map $\mu\in \cM(\P(A),S)$ will be 
called logarithmic if 
$$
\mu(a\cdot a')=\mu(a)+\mu(a'),
$$
for all $a,a'\in A\setminus 0$.
The set of such maps will be denoted by $\LF(\bP(A),S)$.
\end{defn}

\begin{rem}
\label{rem:topo}
In our applications, $A$ and $S$ are endowed with topologies. 
We will always consider {\em continuous} maps, so that 
the notations $\cM, \AF, \LF$ etc. 
stand for spaces of continuous maps, 
subject to the above properties. 
\end{rem}

\begin{defn}
\label{c:pair}
Let $S$ be a ring. Two maps 
$\mu,\mu'\in \cM(\bP(A),S)$ will be called a $c$-pair (commuting pair) if  
\begin{itemize}
\item  $\mu,\mu'\in  \LF(\bP(A),S)$;
\item for all two-dimensional $\F$-subspaces $B\subset A$
there exist constants $s,s',s''$ (depending on $B$) 
with $(s,s')\neq (0,0)$ such that for all $b\in B$ one has
$$
s\mu_B(b)+s'\mu_B'(b)=s''.
$$ 
\end{itemize}
\end{defn}

\begin{thm}
\label{thm:cp}
Let $A$ be an $\F$-algebra and 
$\mu,\mu'\in \cL(\P(A),\Z_{\ell})$ nonproportional maps
forming a $c$-pair. Then there exist constants 
$c,c'\in \Z_{\ell}$, $(c,c')\neq (0,0)$
such that $c\mu+c'\mu'\in \AF(\bP(A),\Z_{\ell})$.   
\end{thm}

\begin{proof}
This is a special case of the main theorem of \cite{BT}.
We outline its proof, since the result 
is crucial for our applications. 

\

{\em Step 1.} (Lemma 3.1.2 in \cite{BT})  
$\mu\in \AF(\bP(A),S)$ iff for all $h\,:\, S\ra \Z/2\Z$
one has $h\circ \mu\in \AF(\bP(A),\Z/2\Z)$.

\

{\em Step 2.} (Lemma 3.2.1 and Proposition 2.6.1 in \cite{BT})
$\mu\in \AF(\bP(A),S)$ iff for all  $B\subset A$ 
with $\dim(B)\le 2$ the restriction $\mu_B\in \AF(\bP(B),S)$. 

\

{\em Step 3.}  
(Proposition 4.3.1 in \cite{BT})
Assume that $\dim(A)\ge 3$ and that $\mu',\mu''\in \LF(\bP(A),\Z_{\ell})$ 
is a $c$-pair such that $\langle \mu',\mu''\rangle_{\Z_{\ell}}$
does not contain an $\AF$-map. 
Then there is a $B\subset A$ with $\dim(B)=3$ 
such that there exist no (nonzero) 
$\AF$-maps $\mu_B\in \langle \mu'_B,\mu_B''\rangle_{\Z_{\ell}}$. 

\

{\em Step 4.} (Lemma 4.3.2 in \cite{BT})
Let $B=\F^3$ and $\mu',\mu''\in \cM(\bP(B),\Z_{\ell})$ be a $c$-pair.
Then the map
$$
\begin{array}{ccccc}
\varphi & : &  \bP(B)  & \ra     &  \bA^2(\Z_{\ell})\\
        &   &     b    & \mapsto &  (\mu'(b),\mu''(b)),
\end{array}
$$
has the following property: the image of every 
$\bP^1\subset \bP(B)$ is contained in an affine line in $\bA^2(\Z_{\ell})$.   
By affine geometry constructions in Section 4.1 of \cite{BT},  
the image $\varphi(\bP(B))$ is 
contained in a union of an affine line and (possibly) 
one more point. Using this, and
assuming that there are no nonzero
$\AF$-maps $\mu\in \langle \mu',\mu''\rangle_{\Z_{\ell}}$, we conclude that
there exist nonconstant (and nonproportional) maps 
$\phi', \phi''\in  \langle \mu',\mu'',1\rangle_{\Z_{\ell}}$
and an $h\,:\, \Z_{\ell}\ra \Z/2\Z$
such that {\em neither} of the three maps
$$
\begin{array}{c}
h\circ \phi',\\
h\circ \phi'',\\
h\circ \phi'+h\circ \phi''
\end{array}
$$ 
is an $\AF$-map on $B$ (with values in $\Z/2\Z$).   

\

{\em Step 5.} 
Rename the three maps to 
$\mu',\mu'',\mu'+\mu''\in \cM(\bP(B),\Z/2\Z)$. 
The image $\varphi(\bP(B))$ is contained in 3 points 
($(0,0)$, $(0,1)$ and $(1,0)$).
Consequently, $\bP(B)$ contains lines of (at most) three types
($(00),(01),(10)$). 
We may assume that there are at least two lines of each type and each line 
contains at least two points of each type (otherwise, it is easy to find 
an $\AF$-map among  $\mu',\mu'',\mu'+\mu''$, leading to a contradiction). 
We may also assume that 
$\mu'$, resp. $\mu''$, $\mu'+\mu''$, is constant on lines of type
$(00)$, resp. $(01)$, $(10)$.   
The projective geometry constructions in Section 4.2 of \cite{BT}
show that one of the three maps has the following property:
on every line $ \bP(C)\subset \bP(B)$ (where it is nonconstant) it
satisfies the functional equation of Step 6. 

\

{\em Step 6.} (Lemma 2.5.6 in \cite{BT})
Let $C=\F^2$ and $\mu\in \cM(\bP(C),\Z/2\Z)$ be such that
$C$ has a basis $(c,b)$ with 
$\mu(c)=\mu(c+\kappa b)\neq \mu(b)$, for all $\kappa \in \F$.
Then $\mu\in \AF(\bP(C),\Z/2\Z)$. 

Indeed, since $\mu$ is homogeneous, 
the functional equation implies that  
$$
\mu(\kappa c+\kappa' b)=\mu(c)
$$ 
for all $\kappa, \kappa'\in \F^*$. 
Thus $\mu$ is constant on $C\setminus \F b$.

\

Thus at least one of the maps $\mu',\mu'',\mu'+\mu''$
is an $\AF$-map on {\em all} lines in $\bP(B)$, 
hence an $\AF$-map on the whole $\bP(B)$, contradiction to Step 4.
\end{proof}

\section{Galois groups}
\label{sect:gal}

Let $k$ be an algebraic closure of $\F$, $K=k(X)$ the function field of an
algebraic variety  $X$ over $k$, $\G^a_K$ the abelianization of the
pro-$\ell$-quotient $\G_K$ of the Galois group  of a separable
closure of $K$, 
$$
\G^c_K=\G_K/[[\G_K,\G_K],\G_K]\stackrel{\pr}{\lra} \G^a_K
$$
its canonical central extension and $\pr$ the natural projection.
By our assumptions, $\G^a_K$ is a torsion free $\Z_{\ell}$-module.

\begin{defn}
\label{defn:lift}
We say that $h,h'\in \G^a_K$ form a {\em commuting pair}
if for some (and therefore any) 
of their preimages $\tilde{h},\tilde{h}'$ one has
$
[\tilde{h},\tilde{h}']=0.
$
A subgroup $\cH$ of $\G^a$ is called {\em liftable}
if any two elements in $\cH$ form a commuting pair. 
\end{defn}

\begin{defn} 
\label{defn:fan}
A {\em fan} $\Sigma_K=\{ \sigma\}$ on $\G^a_K$ is 
the set of all those topologically noncyclic
liftable subgroups $\sigma\subset\G^a_K$
which are not properly contained in any other 
liftable subgroup of $\G^a_K$. 
\end{defn}

\begin{rem}
Even though the group $\G^c_K$ depends on the choice of
a separable closure of $K$, the fan $\Sigma_K$ does not. 
\end{rem}

\begin{rem}
\label{rem:2}
For function fields $E/k$ of curves
there are no topologically noncyclic liftable subgroups of $\G^a_E$ and 
$\Sigma_E=\emptyset$.  
For function fields $K/k$ of surfaces all groups $\sigma\in \Sigma_K$
are isomorphic to (torsion free) $\Z_{\ell}$-modules of rank 2 
(see Section~\ref{sect:af-val}).  
\end{rem}

\begin{nota}
\label{nota:kkk}
Let 
$$
\boldsymbol{\mu}_{\ell^n}:=\{ \sqrt[\ell^n]{1}\,\}
$$
and 
$$
\Z_{\ell}(1) =\lim_{n\ra \infty} \boldsymbol{\mu}_{\ell^n}
$$
be the Tate twist of $\Z_{\ell}$.
Write 
$$
\hat{K}^*:=\lim_{n\ra \infty} K^*/ (K^*)^{\ell^n}
$$
for the multiplicative group of (formal) rational functions on $X$.
\end{nota}

\begin{thm}[Kummer theory]
\label{thm:ga} 
For every $n\in \N$ we have a pairing 
$$
\begin{array}{ccl}
\G^a_K/\ell^n\times K^*/(K^*)^{\ell^n} & \ra & \boldsymbol{\mu}_{\ell^n} \\
   (\mu,f)                          & \mapsto & [\mu,f]_n:= \mu(f)/f
\end{array}
$$
which extends to a nondegenerate pairing 
$$
[ \cdot ,\cdot ]\,:\, 
\G^a_K\times \hat{K}^*\ra \Z_{\ell}(1).
$$
\end{thm}

\begin{rem}
\label{erma:choice}
Since $k$ is algebraically closed and $\ell\neq  p$ we can 
choose a noncanonical isomorphism 
$$
\Z_\ell\simeq \Z_{\ell}(1).
$$ 
From now on we will fix such a choice.
\end{rem}

\section{Valuations}
\label{sect:val}

In this section we recall basic facts concerning valuations 
and valued fields (we follow \cite{bourb}).

\begin{nota}
\label{nota:vg}
A {\em value group}
is a totally ordered (torsion free) abelian group. 
We will denote by $\Ga$ a value group, 
use the additive notation $''+''$ 
for the group law and $\ge$ for the order.
We have 
$$
\Ga=\Ga^{+}\cup \Ga^{-}, \,\,\, \Ga^{+}\cap \Ga^{-}=\{0\}\,\,\, 
{\rm and} \,\,\, \gamma 
\ge \gamma ' \,\,\, {\rm iff} \,\,\,\gamma-\gamma'\in \Ga^{+}.
$$  
We put 
$\Ga_{\infty}=\Ga\cup \{ \infty\}$ and make it to a totally ordered group 
through the conventions
$$
\gamma<\infty, \,\, 
\gamma+\infty=\infty+\infty=\infty, \,\,\,  \forall \gamma\in \Ga. 
$$
\end{nota}

\begin{exam}
\label{exam:lexi}
A standard value group is $\Z^n$ with the {\em lexicographic} ordering.
\end{exam}

\begin{defn}
\label{defn:valu}
A (nonarchimedean) {\em valuation} $\nu=(\nu,\Ga_{\nu})$ on a field $K$
is a pair consisting of a  totally ordered 
abelian group $\Ga_{\nu}$ (the {value} group) 
and a map
$$
K\ra \Ga_{\nu,\infty}
$$
such that
\begin{itemize}
\item $\nu\,:\, K^*\ra \Ga_{\nu}$ is a surjective homomorphism;
\item $\nu(\kappa+\kappa')\ge 
\min(\nu(\kappa),\nu(\kappa'))$ for all $\kappa,\kappa'\in K$;
\item $\nu(0)=\infty$.
\end{itemize}
\end{defn}

\begin{rem}
\label{rem:zeta}
In particular, since $\Ga_{\nu}$ 
is nontorsion,  
$\nu(\zeta)=0$ for every element $\zeta$ of finite order in $K^*$. 
\end{rem}

A valuation is called {\em trivial} if $\Ga=\{0\}$. 
In our applications, 
$K=k(X)$, where $k$ is an algebraic closure of
the finite field $\F$ and 
$X$ is an algebraic variety over $k$ (a {\em model} for $K$). 
Since every element in $k$ is torsion, 
every valuation of $K$ restricts to a trivial valuation on $k$.

\begin{nota}
\label{nota:ko}
We denote by $K_{\nu}$, $\mo_{\nu}, \mm_{\nu}$ and 
$\KK_{\nu}$ 
the completion of $K$ with respect to $\nu$,
the ring of $\nu$-integers in $K$,
the maximal ideal of $\mo_{\nu}$ and the residue field
$$
\KK_{\nu}:=\mo_{\nu}/\mm_{\nu}.
$$
If $X$ (over $k$) is a model for $K$ then
the {\em center} $\mc(\nu)$ of 
a valuation is the irreducible subvariety 
defined by the prime ideal $\mm_{\nu}\cap k[X]$ 
(provided $\nu$ is nonnegative on $k[X]$). 
\end{nota}

It is useful to keep in mind the following exact sequences:
\begin{equation}
\label{eqn:1}
1\ra \mo_{\nu}^*\ra K^*\ra \Ga_{\nu}\ra 1
\end{equation}
and
\begin{equation}
\label{eqn:2}
1\ra (1+\mm_{\nu})\ra \mo_{\nu}^*\ra \KK_{\nu}^*\ra 1.
\end{equation}

\begin{nota}
\label{nota:ine}
We denote by $\I^a_{\nu}\subset \D^a_{\nu}\subset \G^a_{K}$ 
the images of the inertia and the decomposition group of
the valuation $\nu$ in $\G^a_K$.
There is a (canonical) isomorphism 
$$
\D^a_{\nu}/\I^a_{\nu}\simeq \G^a_{\KK_{\nu}}.
$$
\end{nota}

\begin{defn}
\label{defn:div}
Let $K$ be the function field 
of an algebraic variety over $k$.
A valuation $\nu$ of $K$ is {\em positive-dimensional} if 
its residue field $\KK_{\nu}$ is 
the function field of an 
algebraic variety of dimension $\ge 1$.
It is {\em divisorial} if
$$
\trdeg_{k} \KK_{\nu}=\dim(X)-1.
$$
\end{defn}

\begin{nota}
\label{nota:cv}
We let $\Val_K$ be the set 
of all nontrivial (nonarchimedean) valuations of $K$ and 
$\DVal_K$ the subset of its divisorial valuations. If  
$\nu\in \DVal_K$ is realized by a divisor on a model $X$ of $K$ 
(see Example~\ref{exam:suff})
we sometimes write $\I^a_D$, resp. $\D^a_D$, for the corresponding 
inertia, resp. decomposition group.  
\end{nota}

\begin{exam}
\label{exam:curve}
Let $C$ be a smooth curve over $k$ and $E=k(C)$. 
Every point $q\in C$ defines a nontrivial valuation $\nu_q$ on $E$
(the order of a function $f\in E^*$ at $q$). 
Conversely, every nontrivial valuation $\nu$ 
on $E$ defines a point $q:=\mathfrak c(\nu)$ on $C$. 
\end{exam}

\begin{exam}
\label{exam:suff}
Let $X$ be a smooth surface over $k$ and $K=k(X)$. 
\begin{itemize}
\item 
Every curve $C\subset X$ defines a valuation $\nu_C$ on 
$K$ with value group $\Z$. Conversely, every 
valuation $\nu$ on $K$ with value group isomorphic to $\Z$ defines
a curve on {\em some} model $X'$ of $K$.
\item
Every flag $(C,q)$, (curve, point on this curve), 
defines a valuation $\nu_{C,q}$ on 
$K$ with value group $\Z^2$. 
\item There exist (analytic) valuations on $K$ 
with value group a subgroup of $\Z^2$ and center supported in a point
(on every model).  
\end{itemize}
\end{exam}

\begin{rem}
\label{rem:surffa}
Every (nontrivial) positive-dimensional valuation 
on the function field of a surface $X/k$ is divisorial. 
\end{rem}

\begin{defn}
\label{defn:comp}
We say that two valuations $(\nu_1,\Ga_1)$ and $(\nu_2,\Ga_2)$ 
are {\em compatible}
if there exists a valuation $(\nu,\Ga)$ and two
surjective order-preserving homomorphisms of value groups
$\pr_j\,:\, \Ga\ra \Ga_j$ ($j=1,2)$ such that 
for all $\kappa\in K^*$ 
$$
\nu_j(\kappa)=\pr_j(\nu(\kappa)).
$$
Two valuations $(\nu_1,\Ga_1)$ and $(\nu_2,\Ga_2)$ 
are {\em equivalent} if there exists an order preserving isomorphism of 
value groups $\Ga_1\ra \Ga_2$ 
commuting with the homomorphisms
$\nu_1,\nu_2$. 
\end{defn}

We will not distinguish equivalent valuations.

\begin{lemm} 
\label{lemm:compa}
Let $K$ be any field and 
$\nu',\nu''$ two valuations on $K$. 
Then either
\begin{itemize}  
\item $\nu'$ and $\nu''$ 
are compatible and there exists
a valuation $\nu$ with maximal ideal 
$\mm_{\nu} = \mm_{\nu'} + \mm_{\nu''}$
or
\item $\nu'$, $\nu''$ are incompatible,
$$
K = \mm_{\nu'} + \mm_{\nu''}\,\,\, {\rm and} \,\,\, 
K^* = (1 +\mm_{\nu'})\cdot \mo_{\nu''}^*.
$$
\end{itemize}
\end{lemm}

\begin{proof}
See \cite{ZS}.
\end{proof}

\section{A dictionary}
\label{sect:dict}

Write
$$
\begin{array}{rcl}
\cL_K & := & \cL(\P(K),\Z_{\ell}(1)),\\
\AF_K & := & \AF(\P(K),\Z_{\ell}(1)).
\end{array}
$$

\begin{prop}
\label{prop:decomp}
One has the following identifications:
$$
\begin{array}{rcl}
\G^a_K     & = & \cL_K, \\
\D^a_{\nu} & = & 
\{ \mu\in \cL_K\,|\, \mu\,\,\,{\rm trivial }\,\, \, {\rm on}\,\,\,
(1+\mm_{\nu})\},\\
\I^a_{\nu}& = & 
\{ \mu\in \cL_K\,|\, \mu\,\,\,{\rm trivial }\,\, \, {\rm on}\,\,\,
\mo_{\nu}^*\}.
\end{array}
$$
If two nonproportional 
$\mu,\mu'\in \G^a_K$ form a commuting pair then the corresponding maps
$\mu,\mu'\in \cL_K$ form a $c$-pair 
(in the sense of Definition~\ref{c:pair}). 
\end{prop}

\begin{proof}
The first identification is a consequence of Kummer theory \ref{thm:ga}. 
For the second and third identification we use the sequences
\eqref{eqn:1} and \eqref{eqn:2}. 
For the last statement, assume that $\mu,\mu'\in \cL_K$
don't form a $c$-pair. Then there is an $x\in K$ such that
the restrictions of $\mu,\mu'\in\cL_K$ to the subgroup 
$\langle 1,x\rangle $ are linearly independent. Therefore, 
$\mu,\mu'\in \G^a_K$ define a rank 2 liftable subgroup in 
$\G^a_{k(x)}$. Such subgroups don't exist since 
$\G_{k(x)}$ is a free pro-$\ell$-group.
\end{proof}

\begin{exam}
\label{exam:cp}
If $\mu\in \D^a_{\nu}$ and $\al\in \I^a_{\nu}$ then 
$\mu$ and $\al$ form a commuting pair. 
\end{exam}

\begin{prop}
\label{prop:may}
Let $K$ be a field and $\alpha \in \AF_K\cap\cL_K$.
Then there exists a unique valuation $\nu=(\nu_{\al},\Ga_{\nu_{\al}})$ 
(up to equivalence) and a homomorphism $\pr\,:\, \Ga_{\nu_{\al}}\ra \Z_{\ell}$ 
such that 
$$
\al(f)=\pr(\nu_{\al}(f))
$$
for all $f\in K^*$. In particular, $\al\in \I^a_{\nu}$
(under the identification of Proposition~\ref{prop:decomp}).  
\end{prop}

\begin{proof}
Assume that $\al(f)\neq \al(f)$ for some $f,f'\in K$ and 
consider the projective line  $\P^1=\P(\F f + \F f')$. 
Since  $\al$ is a flag map, it is constant outside one point on 
this $\P^1$ so that either $\al(f+f')=\al(f)$  or $=\al(f')$. This 
defines a relation: $f'>_{\al} f$ (in the first case) 
and $f >_{\al} f'$  (otherwise). 
If $\al(f)=\al(f')$ and there exists an $f''$ such that 
$\al(f)\neq \al(f'')$ and $f>_{\al} f'' >_{\al} f'$ then
we put $f >_{\al} f'$. Otherwise, we put $f =_{\al} f'$. 

It was proved in  \cite{BT}, Section 2.4, that 
the above definitions are correct and that $>_{\al}$ is indeed an order
which defines a filtration on the additive group $K$ by subgroups
$(K_{\gamma})_{\gamma\in \Ga}$
such that 
\begin{itemize}
\item 
$K=\cup_{\gamma\in \Ga} K_{\gamma}$ 
and 
\item 
$\cap_{\gamma \in \Ga} K_{\gamma}=\emptyset$, 
\end{itemize}
where $\Gamma$ is the set of equivalence classes with respect to $=_{\al}$.  
Since $\al\in \cL_K$ this order is compatible with multiplication in $K^*$, 
so that the map $K\ra \Gamma$ is a valuation
and $\al$ factors as $K^*\ra\Ga\ra\Z_{\ell}\simeq \Z_{\ell}(1)$.
By \eqref{eqn:1}, $\al\in \I^a_{\nu}$.   
\end{proof}

\begin{coro}
\label{coro:lift}
Every (topologically) noncyclic liftable subgroup of 
$\G^a_K$ contains an inertia element of some valuation. 
\end{coro}

\begin{proof}
By Theorem~\ref{thm:cp}, every such liftable 
subgroup contains an $\AF$-map, 
which by Proposition~\ref{prop:may} belongs to some inertia group.  
\end{proof}

\section{Flag maps and valuations}
\label{sect:af-val}

\begin{lemm}
\label{lemm:am}
Let $\al\in \AF_K\cap \cL_K$, $\nu=\nu_{\al}$ 
the associated valuation and 
$\mu\in \cL_K$. Assume that $\al,\mu$ form a c-pair.
Then 
$$
\mu(1+\mm_{\nu})=\mu(1).
$$
In particular, the restriction of $\mu$ to $\mo_{\nu}$ 
is induced from $\KK_{\nu}$. 
\end{lemm}

\begin{proof}
We have
\begin{itemize}
\item[(1)]  
$\al(\kappa)=0$ for all $\kappa\in \mo_{\nu}\setminus \mm_{\nu}$;
\item[(2)] $\al(\kappa+m)=\al(\kappa)$ for all $\kappa$ and $m$ as above;
\item[(3)]  $\mm_{\nu}$ is generated by $m\in \mo_{\nu}$ such that
$\al(m)\neq 0$.
\end{itemize}
If $m\in \mm_{\nu}$ 
is such that $\al(m)\neq 0$ and 
$\kappa\in  \mo_{\nu}\setminus \mm_{\nu}$  
then $\al$ is nonconstant
on the subgroup $A:=\langle\kappa,m\rangle$. 
Then 
$$
\mu(\kappa+m)=\mu(\kappa).
$$
Indeed, if $\mu$ is 
nonconstant on $A$ 
the restriction $\mu_A$ is proportional to $\al_A$ 
(by the $c$-pair property) and $\al$ satisfies (2).
In particular, for such $m$ we have $\mu(1+m)=\mu(1)$.

If $\al(m)=0$ then there exists $m',m''\in \mm_{\nu}$ such that
$m=m'+m''$ and $\al(m')=\al(m'')\neq 0$. 
Indeed, there exists
an $m'\in \mm_{\nu}$ such that $m>m'>1$ and $\al(m')\neq \al(1)=0$.
Since $\al$ takes only two values on 
the subgroup 
$\langle m',m\rangle\subset \mm_{\nu}$ 
we have 
$$
\al(m'')=\al(-m'+m)=\al(m').
$$
Therefore,   
$$
0=\mu(1 + m')+\mu(1+m'')  = 
\mu(1+ m + m'm'').
$$
Put $\kappa=1+m + m'm''$ and observe that $\al(-m'm'')=-2\al(m')\neq 0$.
By the argument above
$$
\mu(\kappa-m'm'')=\mu(\kappa)=\mu(1+m'+m'')=\mu(1+m),
$$
as claimed.
\end{proof}

\begin{lemm}
\label{lemm:ggg}
Assume that $\al,\al'\in \G^a_K\cap \AF_K$ form a $c$-pair.
Then the associated valuations $\nu$ and $\nu'$ are
compatible.
\end{lemm}

\begin{proof}
If $\nu$ and $\nu'$ are incompatible then, by Lemma~\ref{lemm:compa}, 
$$
K^*=\mo_{\nu}^*\cdot (1+\mm_{\mu'}).
$$ 
By Lemma~\ref{lemm:am}, 
$$
\al(1+\mm_{\nu'})=0 \,\, \text{  and } \,\,\al(\mo_{\nu}^*)=0.
$$ 
This implies that $\al$, and similarly, $\al'$ vanish on $K^*$.
\end{proof}

\begin{coro}
\label{coro:22}
Let $K/k$ be the function field of a surface. Then  
for every $\sigma\in \Sigma_K$ one has
$$
\rk_{\Z_{\ell}}\, \sigma =2.
$$
\end{coro}

\begin{proof}
Otherwise, by Lemma~\ref{lemm:am}, we would get a topologically noncyclic liftable
subgroup in $\G^a_{\KK_{\nu}}$ (for some $\nu\in \Val_K$). 
Since $\KK_{\nu}$ is either $k$ or the function field of some curve, 
we get a contradiction.
\end{proof}

\begin{prop}
\label{prop:centr}
Let $K=k(X)$ be the function field of a surface and 
$\sigma,\sigma'\in \Sigma_K$ be two distinct 
maximal liftable (topologically noncyclic) subgroups of $\G^a_K$
such that $\sigma\cap \sigma'\neq 0$. Then there exists a 
unique divisorial valuation $\nu$ of $K$ such that
\begin{itemize}
\item $\sigma\cap \sigma'=\I^a_{\nu}$;
\item both $\sigma$ and $\sigma'$ are contained in 
$\D^a_{\nu}$.
\end{itemize} 
Moreover, for every divisorial 
valuation $\nu$ of $K$ there
exist $\sigma,\sigma'\in \Sigma_K$ as above. 
\end{prop}

\begin{proof}
First of all, every $\al\in \sigma\cap \sigma'$ is an $\AF$-map
and therefore an inertial element for some valuation $\nu$.  
It commutes only with the associated $\D^a_{\nu}$. 
Since in our case every (nontrivial)
valuation $\nu$ is geometric, its center $\mathfrak c(\nu)$ is
a subvariety on some model $X$ of $K$.
If $\mathfrak c(\nu)$ is a point then 
the rank of $\D^a_{\nu}$ is at most 2.  
However, $\sigma\cup \sigma\subset \D^a_{\nu}$ and has rank  3, 
contradiction. 
Hence $\dim(\mathfrak c(\nu))=1$ and $\nu$ is divisorial.
\end{proof}

\section{Galois groups of curves}
\label{sect:valc}

Throughout this section 
$E/k$ is a 1-dimensional field and $C$ a smooth curve of genus $\mathsf g$ 
with $k(C)=E$.
We have an exact sequence
$$
0\ra E^*/k^*\ra \Div(C)\ra \Pic(C)\ra 0
$$
(where $\Div(C)$ can be identified with 
the free abelian group generated by points in $C(k)$). 
This gives a dual sequence
\begin{equation}
\label{eqn:crrr}
0\ra \Z_{\ell}(\Delta) \ra \cM(C(k),\Z_{\ell})\ra 
\G^a_E\ra \Z_{\ell}^{2\mathsf g}\ra 0,
\end{equation}
with the identifications 
\begin{itemize}
\item $\Hom(\Pic(C),\Z_{\ell})=\Z_{\ell}$ (since $\Pic^0(C)$ is torsion);
\item $\cM(C(k),\Z_{\ell})=\Hom(\Div(C),\Z_{\ell})$ 
is the $\Z_{\ell}$-linear space of maps from $C(k)\ra \Z_{\ell}$;
\item $\Z_{\ell}^{2\mathsf g}=\Ext^1(\Pic^0(C),\Z_{\ell})$.
\end{itemize}
Using this model and the results in Section~\ref{sect:gal}, 
we can interpret 
\begin{equation}
\label{eqn:cm}
\G^a_{E}\subset \cM(C(k),\Q_{\ell})/\text{constant maps} 
\end{equation}
as the $\Z_{\ell}$-linear subspace of 
all maps $\mu\,:\, C(k)\ra \Q_{\ell}$ (modulo constant maps)
such that
$$
[\mu,f]\in \Z_{\ell} \,\, \text{ for all }\,\,  f\in E^*/k^*.
$$
Here $[\cdot ,\cdot ]$ is the pairing: 
\begin{equation}
\label{eqn:inter}
\begin{array}{rcl}
\cM(C(k),\Q_{\ell})\times E^*/k^* & \ra     &       \Q_{\ell} \\ 
       (\mu , f )                 & \mapsto &  [\mu,f]:=\sum_{q} \mu(q)f_q,
\end{array}   
\end{equation}
where $\dv(f)=\sum_q f_q q$.
In detail, let $\gamma\in \G^a_E$ be an element of the Galois group. 
By Kummer theory, $\gamma$ is a homomorphism 
$K^*/k^* \ra \Z_{\ell}=\Z_{\ell}(1)$. 
Choose a point $c_0\in C(k)$. 
For every point $c\in C(k)$, there is an $m_c\in \N$ such that
the divisor $m_c(c-c_0)$ is principal. Define a map 
$$
\begin{array}{rccc} 
\mu_{\gamma}\,:\, & C(k) & \ra      &  \Q_{\ell},\\
                  &  c   &  \mapsto & \gamma(m_c(c-c_0))/m_c.
\end{array} 
$$ 
Changing $c_0$ we get maps differing by a constant map.

In this interpretation, an element of an 
inertia subgroup $\I^a_{w}\subset \G^a_E$
corresponds to a ``delta''-map (constant outside the point $q_w$). 
Each $\I^a_w$ has a canonical (topological) generator $\delta_w$, given by
$\delta_w(f) = \nu_w(f)$, for all $f\in E^*/k^*$. The (diagonal) map 
$\Delta\in \cM(C(k),\Q_{\ell})$ from \eqref{eqn:crrr}
is then given by 
$$
\Delta =\sum_{w\in \Val_E} \delta_w = \sum_{q_w\in C(k)} \delta_{q_w}. 
$$

\begin{defn}
\label{defn:sup}
We say that the support of a subgroup 
$\I\subset \G^a_E$ is $\le s$ and write
$$
|\supp(\I)|\le s
$$
if there exist valuations $w_1,...,w_s\in \Val_{E}$
such that 
$$
\I\subset \langle 
\I^a_{w_1},...,\I^a_{w_s}\rangle_{\Z_{\ell}} \subset \G^a_{E}.
$$ 
Otherwise, we write $|\supp(\I)|> s $.
\end{defn}

\begin{lemm}
\label{lemm:cu}
Let $\I\subset \G^a_{E}$
be a topologically cyclic subgroup such that 
$|\supp(\I)| > s\ge 2$.
Then there exist a finite set 
$\{ f_j\}_{j\in J}\subset E^*$ and an $m\in \N$ such that the map
$$
\begin{array}{cccc}
\psi: &  \G^a_{E}      & \ra        & V:=\oplus_{j\in J} \Z/\ell^m\\
        &      \mu  & \mapsto &  ([\mu, f_j]_{m})_{j\in J}
\end{array}              
$$
has the following property: for every  set 
$\{ w_1,...,w_s\}\subset \Val_E$  
$$
\psi(\I)\not\subset 
\langle \psi(\I^a_{w_1}),...,\psi(\I^a_{w_s})\rangle_{\Z_{\ell}}.
$$
\end{lemm}

\begin{proof}
Let $\iota\in \G^a_E\subset \cM(C(k),\Q_\ell)$ 
be a {\em representative}, as in 
\eqref{eqn:cm},
of a topological generator of $\I$, where $\supp(\I) > s$.
There are three possibilities:
\begin{itemize} 
\item[(1)] 
$\iota(C(k))\subset \Q_{\ell}$ is infinite;
\item[(2)] 
there is a $b\in \iota(C(k))\subset \Q_{\ell}$ such that 
$\iota^{-1}(b)$ is infinite
{\em and} 
there exist at least $s+1$ distinct points 
$q_{s+2},\ldots, q_{2s+2}\in C(k)$
such that $\iota(q_j)\neq b$ for all $j=s+2,\ldots, 2s+2$;
\item[(3)] otherwise: $\iota(C(k))$ is finite, there is a $b$ with $\iota^{-1}(b)$
infinite and there are at most $s$ distinct points with values differing from 
$b$. 
\end{itemize}
In Case (3), $|\supp(\I)|\le s $.
 
In Case (1), choose any set $Q=\{ q_1,...,q_{2s+2}\}\subset C(k)$ 
of points with pairwise distinct values.
In Case (2) choose distinct $q_1,...,q_{s+1}\in \iota^{-1}(b)$ and 
put $Q:=\{ q_1,...,q_{2s+2}\}$.
In both cases, if  $Q'\subset Q$ is any subset of cardinality $|Q'|=s$ then 
$\iota$ is {\em nonconstant} on $Q\setminus Q'$. In particular,  there exist
points $q_{s_1},q_{s_2}  \in Q\setminus Q'$ such that 
\begin{equation}
\label{eqn:iota}
\iota(q_{s_1})\neq \iota(q_{s_2}).
\end{equation} 

We may assume that $\iota(Q)\subset \Z_{\ell}$ 
(replacing $\iota$ by a sufficiently high multiple, if necessary). 
Now we choose an $m''\in \N$ such that all values of $\iota$ on $Q$ remain
pairwise distinct modulo $\Z/\ell^{m''}$.  
Let $\Div^0_Q(C)$ be the abelian group 
of degree zero divisors on $C$ supported in $Q$. 
By Lemma~\ref{lemm:C}, there is an 
$n=n_Q\in \N$ such that $nD$ is principal for every $D\in\Div^0_Q(C)$.
In particular, for every $q_{s_1},q_{s_2}\in  Q$ 
there is a function $f\in E^*$ 
such that $\dv(f) = n(q_{s_1}-q_{s_2})$. 
Write $n=\ell^{m'}\bar{n}$, with  $\gcd(\bar{n},\ell)=1$,
and put $m=m'+m''$.

We have a pairing (Kummer theory)
$$
\begin{array}{rcl}
\G^a_E\times n\Div^0_Q(C)& \ra &  \Z/\ell^{m}\\
         (\mu,f)      & \mapsto &  [\mu,f]_{m}.
\end{array}
$$
Notice that $[\I^a_{w}, f]=0$ for all 
$w$ with $q_w\notin Q$ and all $f\in E^*$ supported in $Q$. 
Further, for every $Q'\subset Q$ with $|Q'|=s$ and points
$q_{s_1},q_{s_2}\in Q\setminus Q'$ as in \eqref{eqn:iota} 
there is an $f\in E^*$ with divisor
$\dv(f)=n(q_{s_1}-q_{s_2})$ 
such that 
$$
[\iota,f]=n\cdot (\iota(q_{s_1})-\iota(q_{s_2})) \neq 0 \mod \ell^{m}
$$ 
and 
$$
[\I^a_{w'},f] =0
$$
for all $\I^a_{w'}$ of $q'\in Q'$.
Let $\{ f_j\}_{j\in J}$ be a basis for $\ell^m\cdot \Div^0_Q(C)$, 
with $f_j\in E^*$. The map
$$
\begin{array}{rcl}
\psi\,:\, \G^a_E & \ra &  \oplus_{j\in J} \Z/\ell^{m} \\
             \mu & \mapsto  & ([\mu,f_j]_{m})_{j\in J}
\end{array}
$$
satisfies the required properties.
\end{proof}

The next step is an {\em intrinsic} definition of 
inertia subgroups 
$$
\I^a_w\subset \D^a_{\nu}/\I^a_{\nu}=\G^a_{k(C)}.
$$
We have a projection 
$$
\pi_{\nu}\,:\, \G^a_K\ra \G^a_{K}/\I^a_{\nu}
$$ 
and an inclusion
$$
\G^a_{\KK_{\nu}}=\D^a_{\nu}/\I^a_{\nu}\hookrightarrow \G^a_{K}/\I^a_{\nu}
$$

\begin{prop}
\label{prop:intri} 
Let $\nu$ be a divisorial valuation of $K$. 
A subgroup 
$$
\I\subset \D^a_{\nu}/\I^a_{\nu}
$$
is the inertia subgroup
of a divisorial valuation of $k(C)=\KK_{\nu}$ iff
for every homomorphism 
$$
\psi\,:\, \G^a_{K}/\I^a_{\nu}\ra V
$$
onto a finite abelian group $V$   
there exists a divisorial valuation $\nu_\psi$ such that
$$
\psi(\I)=\psi \circ \pi_{\nu} (\I^a_{\nu_{\psi}}).
$$
\end{prop}

\begin{proof}
Let $C$ be the smooth model for $\KK_{\nu}=k(C)$, 
$$
\I=\I^a_w\subset \D^a_{\nu}/\I^a_{\nu}
$$ 
the inertia subgroup of a divisorial valuation of $k(C)$
corresponding to a point $q=q_w\in C(k)$ and  
$$
\psi\,:\, \G^a_{K}/\I^a_{\nu}\ra V
$$
a homomorphism onto a finite abelian group. 
Since $\G^a_{K}$ is a pro-$\ell$-group,
we may assume that 
$$
V=\oplus_{j\in J} \Z/\ell^{n_j},
$$ 
for some $n_j\in \N$. Let $n=\max_j(n_j)$. 
By Kummer theory, 
$$
\Hom(\G^a_K,\Z/\ell^n)=K^*/(K^*)^{\ell^n}
$$ 
so that $\psi$ determines elements 
$$
\bar{f}_j\in K^*/(K^*)^{\ell^{n}}
$$
(for all $j\in J$). 
Choose functions $f_j$ projecting to $\bar{f}_j$.
They define a finite set of divisors 
$D_{ij}$ on $X$.  Moreover, $f_j$ 
are not simultaneously constant on $C$
(otherwise, $\psi(\G^a_{k(C)})=\psi(\I^a_{k(C)})$).  
Changing the model $\tilde{X}\ra X$, if necessary,
we may assume that  
\begin{itemize}
\item $C$ is smooth (and irreducible);
\item there exists exactly one irreducible 
component $D$ in the full preimage
of $\cup D_{ij}$ which intersects $C$ in $q$. 
Moreover, this intersection is transversal  
\end{itemize}
(see Section~\ref{sect:basicalg}). 
Then the image of $\I^a_{D}$ under $\psi$ 
is equal to the image of $\I^a_w$.

\

Conversely, we need to show that if $\I\neq \I^a_w$ 
(for some $w\in \DVal_{\KK_{\nu}}$), then 
there exists a homomorphism 
$$
\psi\,:\, \G^a_K/\I^a_{\nu}\ra V
$$
onto a finite abelian group $V$ such that for all $\nu'\in \DVal_K$ 
one has
$$
\psi(\I)\neq \psi\circ \pi_{\nu} (\I^a_{\nu'}).
$$

We consider two cases
\begin{itemize}
\item[(1)] there exist two points $q,q'\in C(k)$ such that
$ \I\subset     \langle \I^a_w,\I^a_{w'}\rangle$;
\item[(2)] otherwise.
\end{itemize}

{\em Case 1.} There exists a rational map 
$\pi\,:\, X\ra \P^1$ such that
its restriction 
$$
\pi\,:\, C\ra \P^1
$$ 
is surjective, unramified
at $q,q'$ and $\pi(q) \neq \pi(q')$.
Under the induced map of Galois groups  
$$
\pi_*(\I)\subset \langle \I^a_{\pi(w)},\I^a_{\pi(w')}\rangle 
$$ 
but is not contained in either $\I^a_{\pi(q)}$ or $\I^a_{\pi(q')}$.
Thus there exist a finite abelian group $V$ and a map 
$\psi \,:\,\G^a_{k(\P^1)} \ra V$ such that
$\psi(\I)\notin \psi(\I^a_{w''})$ for any $q''\in \P^1$.
It follows that
$$
\psi\circ \pi_*(\I)\notin \psi\circ \pi_*(\I^a_{\nu})
$$ 
for any $\nu\in \DVal_K$.

\

{\em Case 2. }
By Lemma~\ref{lemm:cu}, there exist 
a finite set of functions $\bar{f}_j\in k(C)$, 
with support in a finite set $Q=\{ q_0,...,q_s\}\subset C(k)$, and an 
$m\in \N$ such that the homomorphism
$$
\begin{array}{rcl}
\bar{\psi}\,:\, \G^a_{k(C)}& \ra &  V=\oplus_{j\in J}\Z/\ell^m\\
                    \mu  & \mapsto & ([\mu,\bar{f}_j]_m)_{j\in J}
\end{array}
$$
has the property that for all $w,w'\in \DVal_{k(C)}$
$$
\psi(\I)\not\in \langle \psi(\I^a_w),\psi(\I^a_{w'})\rangle_{\Z_{\ell}}.
$$
Next we choose a model for $X$ and $C$ as in Lemma~\ref{lemm:cc}.
In particular, there exist functions  $g_j$ with divisor
$$
\dv(g_j)=n\cdot (D_j-D_0)+(H_j-H_j')
$$
such that all the 
divisors are irreducible, with transversal intersections
and $\dv(g_j)|C=n(q_j-q_0)$.  
These functions $g_j$ define 
a homomorphism
$$
\psi\,:\, \G^a_K/\I^a_{\nu}\ra V.
$$
If $D$ is a divisor on $X$ then $\psi\circ \pi_\nu(\I^a_D)=0$ unless
$D = D_{j}$ for some $j$. In this
case $\psi\circ \pi_\nu(\I^a_{D_j}) = \psi(\I^a_{w_j})$.

Let $\nu'\in \DVal_K$ and $\mathfrak c(\nu')\subset X$ be its center on $X$. 
There are three cases:
\begin{itemize}
\item $\mathfrak{c}(\nu')\not\subset D_j$ for any $j$:
then $\psi\circ \pi_{\nu}(\I^a_{\nu'}) = 0$;
\item  $\mathfrak{c}(\nu')\in D_j^0$, where 
$D_j^0=D_j\setminus (\cup_{j'\neq j} D_j\cap D_{j'})$: 
then  
$$
\psi\circ \pi_{\nu}(\I^a_{\nu'}) \subset\psi(\I^a_{w_j});
$$
\item   $\mathfrak{c}(\nu')\in D_j\cap D_{j'}$  for some $j,j'$: 
then  
$$
\psi\circ \pi_{\nu}(\I^a_{\nu'}) \subset 
\langle \psi(\I^a_{w_j}),
\psi(\I^a_{w_{j'}})\rangle_{\Z_{\ell}}.
$$
\end{itemize}
All three possibilities contradict our assumptions.
\end{proof}

\begin{lemm}
\label{lemm:gen-type}
Let $C/k$ be a curve and $E=k(C)$ its function field. 
Then ${\mathsf g}(C)\ge 1$ iff there exists a homomorphism 
from $\G^a_E$ to a finite (abelian) 
group which maps all inertia elements to $0$.
\end{lemm}

\begin{proof}
Indeed, every curve of genus $\ge 1$ over a finite field of 
characteristic $p$ has unramified coverings of degree 
$\ell$. These coverings define maps of Galois groups,  
which are trivial on all inertia elements. 
If $C$ is rational then $\G^a_{E}$, and hence its image under
every homomorphism (onto any finite group), 
is generated by inertia elements (see the 
exact sequence \eqref{eqn:crrr}). 
\end{proof}

\begin{rem}
\label{rem:gen-tt}
Combining this with Proposition~\ref{prop:intri} we can 
decide in purely Galois-theoretic terms which
divisorial valuations of $K$ correspond to 
nonrational (irreducible) curves $C$
on some model $X$ of $K$. We call such valuations {\em nonrational}.  
\end{rem}

\section{Valuations on surfaces}
\label{sect:recc}

Let $X$ be a smooth surface over $k$, $K=k(X)$ its function field
and $\nu$ a divisorial valuation of $K$.
We have a well-defined (bilinear, with respect to multiplication)
residue map
\begin{equation}
\label{eqn:ress}
\begin{array}{ccc}
K^*\times K^* & \ra      & \KK_{\nu}/k^*\\
  f,g         & \mapsto  & f^{\nu(g)}/g^{\nu(f)}. 
\end{array}
\end{equation}
On a smooth model $X$ of $K$, where $\nu=\nu_D$
for some divisor $D\subset X$, we can define 
\begin{equation}
\label{eqn:res-1}
\Res_{\nu}=\Res_D \,:\, K^*\times K^* \ra  \KK_{\nu}/k^*
\end{equation}
as follows:
\begin{itemize}
\item $\Res_\nu(f,g) = 0$ if both $f,g$ are invertible on $D$;
\item $\Res_\nu(f,g) = f_D^{m}$
if $f$ is invertible ($f_D$ is the restriction to $D$) 
and $g$ has multiplicity $m$ along $D$;
\item $\Res_\nu(f,g) = (f^{m_g}/{g}^{m_f})_D$ 
in the general case, when 
$f,g$ have multiplicities $m_f,m_g$, respectively. 
\end{itemize}
The definition does not depend on the choice of the model.

\begin{lemm}
\label{lemm:fg}
For $f,g \in K^*$ 
$$
\Res_{\nu} (f,g) =0 \,\,\, \forall \nu\in \DVal_K \,\Longleftrightarrow
f,g\in E=k(C)\subset K \text{ for some  curve } C .
$$ 
\end{lemm}
 
\begin{proof}
($\Leftarrow$) On an appropriate model $X$ we have 
$\nu=\nu_D$ for a divisor $D\subset X$ and 
$\pi\,:\, X\ra C$ is regular and flat with irreducible 
generic fiber (and $f,g\in k(C)^*$). 
By definition, $\Res_\nu(f,g) = 0$ if $D$ is not in the fiber
of $\pi$. If $D$ is in the fiber then
there is a $t\in k(C)^*, \nu_D(t)\neq 0 $ such that
both $ft^{m_f}, gt^{m_g}$ are regular and constant on $D$ 
(for some $m_f,m_g\in \N$) so that $\Res_\nu(f,g) = 0$.
   
($\Rightarrow$) Assume that 
$\Res_\nu(f,g) = 0$ for every $\nu\in \DVal_K$. 
Every nonconstant function $f$ defines a unique map 
(with irreducible generic fiber)
$$
\pi_f\,:\, X\ra C_f
$$
which corresponds to the algebraic closure of $k(f)$ in $K$
(we will say that $f$ is induced from $C_f$).
We claim that $\pi_f  = \pi_g$. 

Since $f$ is induced from $C_f$, we have
$$
\dv(f)=\sum_{q\in Q}a_qD_q,
$$
where $Q\subset C_f(k)$ is finite and $D_q=\pi^{-1}(q)$.
Then $D_q^2 = 0$ and $D_q$ is either a multiple of a fiber of $\pi_g$
or it has an irreducible component $D\subset D_q$ which dominates $C_g$  
(under $\pi_g$).
In the second case, $\nu_{D}(f) \neq 0$, while $\nu_D(g)=0$ and
$g$ is nonconstant on $D$. Hence $\Res_{D}(f,g)\neq 0$, contradiction.
Therefore, all $D_q$ are (multiples of) fibers of $\pi_g$
and $f$ is induced from $C_g$. Hence $C_f=C_g$ and $\pi_f=\pi_g$.
\end{proof}

\section{$\ell$-adic analysis: generalities}
\label{sect:ella}

Let $K=k(X)$ be the function field of a smooth algebraic variety $X$ over $k$.
We have an exact sequence 
\begin{equation}
\label{eqn:seqq}
0\ra K^*/k^*\stackrel{\rho_X}{\lra} \Div(X)\stackrel{\pic}{\lra} \Pic(X)\ra 0, 
\end{equation}
where $\Div(X)$ is the group of (Weil or Cartier) divisors of $X$.
We will identify an element $f\in K^*/k^*$ with its image under $\rho_X$.  
Let 
$$
\widehat{\Div}(X):=\{ D=\sum_{m\in M} \hat{a}_m D_m\}, \,\,\, \text{ resp. }\,\,\,
\widehat{\Div}_{\rm nr}(X) \subset\widehat{\Div}(X),
$$
be the group of divisors (resp. nonrational divisors) 
with {\em rapidly decreasing coefficients}:
\begin{itemize}
\item $M$ is a countable set;
\item  for all $r\in \Z$ the set
$$
\{ m\,|\, |\hat{a}_{m}|_{\ell} \le r\}
$$
is finite;
\item for $D\in\widehat{\Div}_{\rm nr}(X)$, all $D_m$ are nonrational.  
\end{itemize}
Clearly, the group of {\em finite} $\ell$-adic divisors 
$$
\Div(X)_{\ell}:=\Div(X)\otimes_{\Z}\Z_{\ell}\subset \widehat{\Div}(X).
$$
Every element 
$$
\hat{f}\in \hat{K}^*= \lim_{n\ra \infty} K^*/(K^*)^{\ell^n}
$$ 
has a  representation
$$
\hat{f}=(f_n)_{n\in \N} 
\,\, {\rm or }\,\, f=f_0f_1^{\ell}f_2^{\ell^2}\cdots ,
$$
with $f_n\in K^*$. 
We have homomorphisms 
$$
\begin{array}{rcl}
\hat{\rho}_{X} \,: \,\hat{K}^* & \ra  & \widehat{\Div}(X),\\
\hat{f}                  & \mapsto &  
\dv(\hat{f}):=\sum_{n\in \N}\ell^n\cdot\dv(f_n)=
\sum_{m}\hat{a}_m D_m, \\
  & & \\
\hat{\rho}_{X,{\rm nr}} \,: \,\hat{K}^* & \ra  
& \widehat{\Div}(X)\stackrel{{\rm pr}}\lra  
\widehat{\Div}_{\rm nr}(X) ,
\end{array}
$$
where $D_m\subset X$ are irreducible divisors, 
$$
\hat{a}_m=\sum_{n\in \N} a_{nm}\ell^n\in \Z_{\ell},
$$
with $a_{nm}\in \Z$, and 
$$
\dv(f_n)=\sum_{m}a_{nm}D_m.
$$
Here $\dv(f_n)$ is the {\em Cartier} divisor of $f_n$ and 
$\sum_{m}a_{nm}D_m$ is its image in the group of {\em Weil} divisors.
Every $\nu\in \DVal_K$ gives rise to a
homomorphism
$$
\nu\,:\, \hat{K}^*\ra \Z_{\ell}
$$
and a residue map
$$
\hat{\Res}_{\nu}\,:\, \hat{K}^*\times \hat{K}^*\ra \hat{\KK}_{\nu}.
$$
On a smooth model $X$, where $\nu=\nu_D$ for some divisor $D\subset X$, 
$\nu(\hat{f})$ is the $\ell$-adic coefficient at $D$ of $\dv(\hat{f})$, 
while 
$\hat{\Res}_{\nu}$ is the natural generalization of 
\eqref{eqn:ress}.
We say that two elements $\hat{f},\hat{g}\in \hat{K}^*$ commute if
%\begin{itemize}
%\item 
$\hat{\Res}_{\nu}(\hat{f},\hat{g})=0$, for all divisorial $\nu$. 
%\item $\supp_K(\hat{f})\cap \supp_K(\hat{g})=\emptyset$.
%\end{itemize}

\begin{nota}
\label{nota:supp}
We put
$$
\begin{array}{ccccl}
\supp_K(\hat{f})& := \{ &  \nu\in \DVal_K & | &  \hat{f} \,\,\, 
{\rm nontrivial\,\,\, on} 
\,\,\,\I^a_{\nu}\,\,\};\\
\supp_X(\hat{f})& := \{ &  D_m    & | &   \hat{a}_m\neq 0\,\,\}.\\
\end{array}
$$ 
\end{nota}

\begin{defn}
\label{rem:alt}
We say that $\hat{f}$ has {\em finite nonrational support} if 
the set of nonrational $\nu\in \supp_K(\hat{f})$ is finite
(see Lemma~\ref{lemm:gen-type} for the definition and 
Galois-theoretic characterization of nonrational valuations).
Let 
$$
\FS(K)\subset \hat{K}^*
$$
be the set of such elements.
\end{defn}

\begin{defn}
\label{defn:fs}
We say that $\hat{f}$ has finite support on the model $X$ if
$\supp_X(\hat{f})$ is finite. Put
$$
\begin{array}{rl}
\FS_X(K) & :=\{ \hat{f} \in \hat{K}^*\,\,|\, 
\rho_X(\hat{f})\in \Div(X)_{\ell}\}.
\end{array}
$$

\end{defn}

\begin{lemm}
\label{lemm:indep}
The definition of $\FS_X(K)$ does not depend on the choice of a
model $X$. 
\end{lemm}

\begin{proof}
For any two models $X',X''$ we can find
a model $X$ dominating both.  The difference 
between  the sets of irreducible 
divisors $\Div(X')$, resp. $\Div(X'')$,
and $\Div(X)$ is finite (and consists only of rational curves). 
\end{proof}

\begin{coro}
\label{coro:almost}
Let $K$ be the function field of a surface $X$ which
contains only finitely many rational curves. Then 
$$
\FS(K)=\FS_X(K).
$$
\end{coro}

In particular, we obtain an intrinsic, Galois-theoretic 
description of $\FS_X(K)$ in this case.
We proceed to give such a description in general.
Note that for $\hat{f}\in \FS(K)$, its nonrational
component $\hat{\rho}_{X,{\rm nr}}(\hat{f})$ 
is independent of the model $X$. More precisely, 
for any birational morphism $X'\ra X$ we can identify 
$\widehat{\Div}_{\rm nr}(X')=\widehat{\Div}_{\rm nr}(X)$. 
Under this 
identification 
$$
\rho_{X',{\rm nr}}(\hat{f}) = \rho_{X,{\rm nr}}(\hat{f}).
$$
Let ${\mathcal F}(K)$ be the set of all  
$f\in K^*/k^*$ such that $\rho_{X,{\rm nr}}(f)\neq 0$ and 
for every rational divisorial valuation $\nu$ and some 
(equivalently, every) model 
$X$ of $K$, where $\nu=\nu_C$ for a rational curve $C\subset X$,
either 
\begin{itemize}
\item
$f_C=1\in k(C)^*/k^*$ or
\item 
$\rho_C(f_C) \neq 0 \mod \ell$.
\end{itemize}

\begin{lemm}
\label{lemm:imp}
The set $\mathcal F(K)$ generates $K^*/k^*$. Moreover, 
for every pair of commuting elements 
$\hat{f},\hat{g}\in \FS(K)$ with disjoint support such that 
there exist $f,g\in {\mathcal F}(K)$ with 
$$
f=\hat{f} \mod (K^*)^{\ell} \,\,\, \text{ and }\,\,\, 
g=\hat{g} \mod (K^*)^{\ell}, 
$$  
one has $\hat{f}\in \FS_X(K)$ and $\hat{g}\in \FS_X(K)$, 
for every model $X$ of $K$.
\end{lemm}

\begin{proof}
Let $y\in K^*$ be a function such that 
the generic fiber of the corresponding map 
$\pi_y\,:\, X\ra  \P^1_y$, 
from some model $X$ of $K$, is an irreducible nonrational curve. 
Such $y$ generate $K^*$. 

For generic quadratic, coprime polynomials $P,Q\in k[y]$,  
the preimage in $X$ of $(0\cup \infty) \subset \P^1$
under the composition of $\pi_y$ with the map 
$$
\begin{array}{cccc}
\phi\,:\, & \P^1_y & \ra     &   \P^1 \\
         &    y   & \mapsto &  f(y):=P(y)/Q(y)
\end{array}
$$
contains at least 4 irreducible smooth fibers of $\pi_y$. 
If $f$ were nonconstant on a rational curve 
$C$ (on some model $X$ of $K$) and $f_C$ were an $\ell$-th power 
then the local ramification indices of $f$ and hence of $y$ 
were divisible by $\ell$. Thus we would have a map 
$\pi_y\,:\, C\ra \P^1_y$ with all local ramification 
indices over 4 points divisible by $\ell$, and by 
Hurwitz' theorem, 
${\mathsf g}(C)>0$, which contradicts the rationality of $C$. 
It follows that $f\in {\mathcal F}(K)$. 
Clearly, such elements $f$ generate $k(y)^*$.

Next, write
$$
\begin{array}{ccc}
\rho_X(\hat{f})& =& 
\sum_{i\in I} n_i D_i   + \ell \sum_{j=1}^{\infty} n_j C_j, \\
\rho_X(\hat{g})&= & 
\sum_{i\in I'} n_i' D_i' + \ell \sum_{j=1}^{\infty} n_j' C_j',
\end{array}
$$
where $I,I'$ are finite sets and  
the second sum is an infinite series over distinct rational
curves $C_j,C_j'\subset X$. 
By assumption, the sets $\{D_i\}_{i\in I}$, 
$\{C_j\}_{j\in \N}$, $\{D_i'\}_{i\in I'}$, 
$\{C_j'\}_{j\in \N}$ are disjoint.

By assumption, $\rho_{\nu_j'}(\hat{f},\hat{g})=0$, 
for all $\nu_j'$ corresponding to $C_j'$.
Since $C_j'$ are rational, this residue equals the residue of 
$f$ on ${C_j'}$,
which is nonzero $\mod \ell$, contradiction.
Thus, if $(\hat{f}, \hat{g}) =0$, 
then $\supp_X(\hat{g})$ is finite and we may put $\hat{g}=g'$.
The restriction of $g'$ to any irreducible component of the divisor
of $\hat{f}$ is identically zero. This implies that $g'$  is a product of
$\ell$-adic powers of elements belonging to the same field $k(y)$ as $f$.
Thus all rational curves in the support of $\hat{f}$
also belong to the fibers of $y$. There are 
finitely many such curves since some fibers contain nonrational curves.
\end{proof}

We have an exact sequence
$$
0\ra \hat{K}^*\stackrel{\hat{\rho}_{X}}{\lra}
\widehat{\Div}(X)\stackrel{\pic_{\ell}}{\lra}
\Pic(X)_{\ell}\ra 0,
$$ 
where $\Pic(X)_{\ell}:=\Pic(X)\otimes \Z_{\ell}$.
We write 
$$
\widehat{\Div}(X)^0\subset \widehat{\Div}(X)
$$
for the group generated by the image $\hat{\rho}_X(\hat{K}^*)$ and identify
an element $\hat{f}\in \hat{K}^*$ with its image.

\begin{lemm}
\label{lemm:DD}
Let $X/k$ be smooth algebraic with $\NS(X)=\Pic(X)$. Let 
$M$ be a {\em finite} set and  
$$
D=\sum_{m\in M} a_mD_m\in \Div(X)_{\ell}:=\Div(X)\otimes_{\Z}\Z_{\ell},\, \,\, 
a_m\in \Z_{\ell}
$$ 
a divisor such that $\pic_{\ell}(D)=0$.  
Then there exist a finite set $I$, functions $f_i\in K^*$ and numbers 
$a_i\in \Z_{\ell}$, linearly independent over $\Z$,  
such that for all $i\in I$
$$
\supp_{X}(f_i)\in \supp_X(D)
$$ 
and 
$$
D = \sum a_i\dv(f_i).
$$ 
\end{lemm}

\begin{proof}
We have a diagram
$$
\begin{array}{cccccc}
\Ker(\pic) & \ra & \oplus_{m\in M} \Z D_m &\stackrel{\pic}{\lra} & 
\Lambda\subset \Pic(X) & \ra 0\\
\downarrow  &  & \downarrow  &   &     \downarrow  &  \\
\Ker(\pic_{\ell}) & \ra & \oplus_{m\in M}\Z_{\ell}D_m &\stackrel{\pic_{\ell}}{\lra} & 
\Lambda_{\ell}\subset \Pic(X)_{\ell} & \ra 0.
\end{array}
$$
Since $\Pic(X)=\NS(X)$ the map $\Pic(X)\ra \Pic(X)_{\ell}$ is injective and  
$$
\rk_{\Z}\La= \rk_{\Z_{\ell}}\La_{\ell} \text{ and }
\rk_{\Z}\Ker(\pic) = \rk_{\Z_{\ell}} \Ker(\pic_{\ell}).
$$
In particular, $\Ker(\pic_{\ell})$ has a basis 
$\{ D_{i}\}_{i\in I}$ (over $\Z_{\ell}$), where 
each $D_i$ is a $\Z$-integral linear combinations of $D_m$ (with $m\in M$)
and is also in $\Ker(\pic)$. It follows that 
$D_i=\dv(f_i)$ for some function $f_i\in K^*$ with support in $D$. 
Finally, we can find a representation 
$$
D=\sum_{i} a_i D_i, 
$$
with  $a_i\in \Z_{\ell}$ linearly independent over $\Z$
(passing to a subset of $I$, if necessary). 
\end{proof}

\section{$\ell$-adic analysis: curves}
\label{sect:curves-ell}

\begin{prop}
\label{prop:cc}
Let $\tilde{k}$ be the closure of a finite field, $\char(\tilde{k})\neq p$,
$C$ a curve over $\tilde{k}$ of genus 
$\mathsf g$ with function field $E=\tilde{k}(C)$ 
and 
$$
\Phi\,:\, \G^a_{k(\P^1)}\ra \G^a_E
$$
be an isomorphism of Galois groups  
inducing an isomorphism on inertia groups 
of divisorial valuation, that is,  
a bijection on the set of such groups and isomorphisms of
corresponding groups.  
Let 
$$
\Phi^* \,:\, \widehat{k(\P^1)}^*\ra \hat{E}^*
$$
be the corresponding dual isomorphism. 
Then $E=\tilde{k}(\P^1)$ and there is a constant $a\in \Z_{\ell}^*$ such that
$\Phi^*(k(\P^1)^*/k^*)=a\cdot E^*/\tilde{k}^*$. 
\end{prop}  

\begin{proof}
Recalling the exact sequence \eqref{eqn:crrr}, we have a commuting diagram

\centerline{
\xymatrix{  
0\ar[r] & \Z_{\ell}\Delta_{C(\tilde{k})} \ar[r] & \cM(C(\tilde{k}),\Z_\ell) \ar[r]
        & \G^a_E \ar[r] \ar[d]& \Z_{\ell}^{2\mathsf g}\ar[r] & 0 \\
0\ar[r] & \Z_{\ell}\Delta_{\P^1(k)} \ar[r] & \cM(\P^1(k), \Z_{\ell}) \ar[r]
        & \G^a_{k(\P^1)} \ar[r]&  0 \ar[r] &  
}
}

\

Since $\Phi$ is an isomorphism on inertia groups $\I^a_w$, for each $w$, 
the {\em sets} $C(\tilde{k})$ and $\P^1(k)$ coincide and 
we get a {\em unique} isomorphism of $\Z_{\ell}$-modules 
$$
\cM(C(\tilde{k}),\Z_\ell)=\cM(\P^1(\tilde{k}), \Z_{\ell}).
$$
In particular, we find that $\mathsf g=0$ and $E=\tilde{k}(\P^1)$. 
Further, we have an induced isomorphism 
$$
\Z_{\ell} (\sum_{w\in \Val_E} \delta_w) =\Z_{\ell} (\sum_{w'\in \Val_{k(\P^1)}} \delta_{w'})
$$
so that
$$
(\sum_{w\in \Val_E} \delta_w)= a (\sum_{w'\in \Val_{k(\P^1)}} \delta_{w'})
$$
for some $a\in \Z_{\ell}^*$. 
This implies that 
$\delta_w=a\delta_{w'}$, for all $w\in \Val_E$ and the corresponding 
$w'\in \Val_{\P^1}$. In particular, for the dual groups we have
$$
E^*/\tilde{k}^* =(K^*/k^*)^{a}, 
$$
where $a\in \Z_{\ell}^*$.  
\end{proof}

\section{$\ell$-adic analysis: surfaces}
\label{sect:surr}

Let $K=k(X)$ be a function field of a smooth surface $X$ over $k$. 
We will need an $\ell$-adic version of Lemma~\ref{lemm:fg} .

\begin{prop}
\label{prop:gff}
Let $\hat{f},\hat{g}\in \FS(K)$ be such that
\begin{itemize}
\item $\Res_\nu(\hat{f},\hat{g})=0$ for every $\nu\in \DVal_K$;
\item $\supp_K(\hat{f})\cap \supp_K(\hat{g})=\emptyset$.
\end{itemize}
Then there is a 1-dimensional field $E=k(C)\subset K$ such that
$\hat{f},\hat{g}\in \hat{E}^*$. 
\end{prop}

\begin{proof}
By Lemma~\ref{lemm:DD}, 
$$
\hat{f} =\prod_{i\in I} f_i^{a_i}, \text{ resp. } 
\hat{g} =\prod_{j\in J} g_j^{b_j},
$$ 
where 
\begin{itemize}
\item $I,J$ are finite sets;
\item $f_i,g_j\in K^*$ for all $i,j$;
\item $a_i\in \Z_{\ell}$ (resp. $b_j\in \Z_{\ell}$) are linearly independent over $\Z$.
\end{itemize}
Fix a valuation $\nu$ and choose a (smooth) model $X$ so that
$\nu=\nu_D$ for some divisor $D\subset X$. 
Then  
$$
\Res_\nu(\hat{f},\hat{g}) = \prod \Res_D(f_i, g_j)
$$
and we can compute it using only those  
pairs $f_i,g_j$ which have $D$ in their support.
In particular, 
$$
\hat{f}^{m_g} / \hat{g}^{m_f} =  \prod (f_i^{a_im_j}/g_j^{ b_j m_i}), 
$$
where $m_j$ (resp. $m_i$) is the order of $g_j$ (resp. $f_i$) on $D$.
This order vanishes unless $D\in \supp(\hat{f})\cup \supp(\hat{g})$.
By assumption, if $D\in \supp(\hat{f})$ then  
$D\not\in \supp(\hat{g})$ (and $n_j =0$) so that
$$
\Res_D(\hat{f},\hat{g})  \in \hat k(D)^*.
$$
Since the nonzero numbers $a_i$ 
are linearly independent over $\Z$ the equality 
$\sum a_im_i = 0$ implies that $m_i = 0$ (for all $i$) and that $g_D \in k^*$.
    
Similarly, $g_D =\prod (g_j)_D^{b_j}$, where
$b_j$ are linearly independent over $\Z$, and 
$g_D \in k^*$ implies that $(g_j)_D \in k^*$ (for all $j\in J$).
It follows that 
$$
\Res_\nu(f_i,g_j) = 0
$$ 
for all $f_i,g_j $ and every valuation  $\nu=\nu_D$.
By Lemma~\ref{lemm:fg}, all $f_i,g_j$ belong to the same 
1-dimensional field $E\subset K$ and hence $\hat{f},\hat{g}\in \hat{E}^*$.    
\end{proof}

\begin{rem}
\label{rem:sati}
For every $f\in K^*$ the element $g=(f+a)(f+b)$
where $a\neq b$ and $ab\neq 0$, satisfies the conditions 
of Proposition~\ref{prop:gff}.  
\end{rem}

\begin{prop}
\label{prop:surfa}
Let $\sK^*\subset \FS(K)\subset \hat{K}^*$ be a subset with the 
following properties:
\begin{itemize}
\item $\sK^*$ is closed under multiplication;
\item $\sK^*\cap \hat{E}^* = a_E \cdot E^*/k^*$
for every 1-dimensional subfield 
$E=k(x)\subset K$, with $a_E \in \Z_{\ell}^*$;
\item there exists a $\nu_0\in \DVal_K$  such that 
$$
\{ [\delta_0,\hat{f}]\,|\, \hat{f}\in \sK^*\}\simeq \Z
$$
for a topological generator $\delta_0$ of $\I^a_{\nu_0}$.
\end{itemize}
Then $\sK^*\subset K^*/k^*\otimes \Z_{(\ell)}$.  
\end{prop}

\begin{proof}
For  $x\in K\setminus k$ let $E=k(x)$ 
be the corresponding 1-dimensional field. 
By assumption, there 
exists an $a_E\in \Z_{\ell}$ such that 
$$
\sK^*\cap \hat{E}^*=a_E\cdot E^*/k^*.
$$ 
If some (any) topological generator $\delta_0$ of $\I^a_{\nu_0}$ is 
not identically zero on $\hat{E}^*$ 
then there exists a (smooth) model $X$, 
where $\nu_0$ is realized by a divisor $D_0$,  
together with a morphism 
$$
X\ra \P^1=\P^1_E
$$
such that $D_0$ dominates $\P^1$. It follows that 
$$
a_E\in \Q\cap \Z_{\ell}^*=\Z_{(\ell)}.
$$  
It remains to observe that every
$x\in K^*$ can be written as a product
$$
x=x'\cdot x''
$$
such that $\delta_0$ is nontrivial on  
both $E'=k(x')$ and $E''=k(x'')$. 
\end{proof}

\begin{coro}
\label{coro:pm}
After a choice of $\delta_0$, for every 1-dimensional 
$E\subset K$ and every $f\in E^*/k^*$
we can Galois-theoretically distinguish its poles from its zeroes. 
\end{coro}

The last essential step is a Galois-theoretic characterization of
the partial projective structure on $\sK^*/k^*$, more precisely,
the characterization of generating elements 
and primary lines in $\sK^*/k^*$ (see Definition~\ref{defn:gener}
and Example~\ref{exam:basic2}).

\begin{lemm}
\label{lemm:unnn}
Let $x\in K^*$ be a generating element, $E:=k(x)$ 
and $r=r(x)\in \N$ the smallest
positive integer such that $x^r\in \sK^*$. 
Then 
\begin{itemize}
\item $r=p^m$ for some $m\in \N$ (with $p=\char(k)$);
\item $(E^*/k^*)\cap (\sK^*/k^*) = (E^{p^m})^*/k^*$;
\item (pointwise) $p^m$-th powers of primary lines
in $E^*/k^*$ coincide with primary lines in $(E^{p^m})^*/k^*$.
\end{itemize}
\end{lemm}

\begin{proof}
The first property follows since
$K/\sK$ is a finite purely inseparable extension, by 
Propositions~\ref{prop:KKK} and \ref{prop:surfa}. 
Next, we claim that a generator $y\in \sK$  
is a $p^m$-th power of a generator of $K$ (for some $m$ depending on $y$).
Indeed, $E:=\ovl{k(y)}^K\subset K$ is a  finite and purely inseparable 
extension of $k(y)$, $E:=k(x)$ (for some $x\in K$). 
Thus 
$$
y = (ax^{p^m} + b)/ (cx^{p^m} + d) = ((a'x+ b')/(c'x+d'))^{p^m}
$$
for some $ m\in \Z$, $a,b,c,d \in k$ and their $p^m$-th roots 
$a',b',c',d' \in k$ (since $k$ is algebraically closed).

In particular, a generator $y\in \sK^*$ is in $E^*\cap \sK^*$ 
(and is the minimal positive power of a generator in $E$ contained in
$E^*\cap \sK^*$). This implies the third property:
the generators of $E^{p^m}$ are $p^m$-th powers
of the generators of $E$. 
\end{proof}
 
\begin{coro}[Definition]
\label{coro:pro-uni}
Assume that $y,y'$ are primitive elements in  
$(E^{p^m})^*\subset \sK^*$ such that
\begin{itemize}
\item $y,y'$ have support in 2 points;
\item the pole of $y$ coincides with the pole of $y'$.
\end{itemize} 
Then (the images of) 
$y,y'$ in $\sK^*/k^*$ are contained in a primary line 
passing through (the images of) $1,y,y'$.
\end{coro}

\begin{proof}
Definition~\ref{defn:sup} and Lemma~\ref{lemm:cu} give a 
Galois-theoretic 
characterization of the notion ``support in 2 points''. 
By Corollary~\ref{coro:pm} we can Galois-theoretically 
distinguish zeroes and poles of $y\in \sK^*/k^*$. 
It remains to apply Lemma~\ref{lemm:unnn}. 
\end{proof}

\

\section{Proof}
\label{sect:proof}

In this section we prove our main theorem:
if 
$$
(\G_{K}^a,\Sigma_K)=(\G_{L}^a,\Sigma_L),
$$ 
where $L$ is a function field over an algebraic closure
of a finite field of characteristic $\neq 2,\ell$, then 
$K$ is a purely inseparable extension of $L$.

\

{\em Step 1.} We have a nondegenerate pairing 
$$
\G_K^a\times \hat{K}^*\ra \Z_{\ell}(1).   
$$
This implies  that $\hat{K}^*=\hat{L}^*$. 

\

{\em Step 2.}
We identify intrinsically 
the inertia and decomposition groups of divisorial valuations:
$$
\I^a_\nu\subset\D^a_\nu\subset \G^a_K:
$$ 
every liftable subgroup
$\sigma\in \Sigma_K$ contains an inertia element
of a divisorial valuation 
(which is also contained in at least one other $\sigma'\in \Sigma_K$). The
corresponding decomposition group is the ``centralizer'' of the 
(topologically) cyclic inertia group 
(the set of all elements which ``commute'' with inertia). 
This identifies $\DVal_K=\DVal_L$.

\

{\em Step 3.} For every  $\nu\in \DVal_K$
we characterize intrinsically  
$$
\I^a_w\subset \D^a_\nu/\I^a_{\nu}
$$
(see Proposition~\ref{prop:intri}). 

\

{\em Step 4.} We distinguish divisorial valuations with
nonrational centers (see Lemma~\ref{lemm:gen-type} and Remark~\ref{rem:gen-tt}).

\

{\em Step 5.} For $\hat{f}\in \hat{K}^*$ we have two notions of support:
$\supp_K(\hat{f})$ (intrinsic) and $\supp_X(\hat{f})$ (depending on a
model $X$) and two notions of finiteness: 
$\hat{f}$ is nontrivial on at most finitely many
nonrational divisorial valuations $\nu$, 
resp. $\hat{f}$ has finite divisorial support on a model.
We defined $\FS(K)\subset \hat{K}^*$ as the set of elements satisfying 
the first notion of finiteness.  
If some (any) model $X$ of $K$ 
contains only finitely many rational curves, both 
notions of finiteness of support 
coincide and one obtains an intrinsic Galois-theoretic
characterization of $K^*/k^* \otimes \Z_{\ell}\subset \hat{K}^*$, 
as elements in $\FS(K)$.  
In general, it may happen that some $g\in L^*/l^*$
has an ``infinite rational tail'' on some (every) 
model $X$ of $K$:  
$$
\rho_X(g)=\rho_{X,{\rm nr}}(g)+\sum_{j\ge 1} n_j C_j,
$$ 
where $C_j$ are irreducible rational curves on $X$. 
In Lemma~\ref{lemm:imp} we show that a many (and consequently, all) 
elements of $L^*/l^*\subset \FS(L)=\FS(K)$  
have finite support on every model $X$ of $K$, and vice versa.
In particular,  $K^*/k^*\otimes \Z_{\ell} = L^*/l^* \otimes \Z_{\ell}$.

\

{\em Step 6.} For every 
pair of elements $\hat{f},\hat{g}\in \FS_X(K)$ satisfying 
\begin{itemize}
\item $\supp_K(\hat{f})\cap \supp_{K}(\hat{g})=\emptyset$;
\item $\varrho_{\nu}(\hat{f},\hat{g})=0$ for all $\nu\in \DVal_K$
\end{itemize}
there exists a subfield $E=k(C)\subset K$ such that
$\hat{f},\hat{g}\in \hat{E}^*$ (Proposition~\ref{prop:gff}).

\

{\em Step 7.} Since $\Pic(X)=\NS(X)$ every such subfield $E=k(x)$ for
some $x\in K^*$. 

\

{\em Step 8.} Proposition~\ref{prop:cc} identifies 
$E^*/k^*$ inside $\hat{E}^*$, up to conformal equivalence.

\

{\em Step 9.} Proposition~\ref{prop:surfa} identifies
$\sK^*:=K^*\cap L^*$ (as a multiplicative group) with a 
multiplicative subgroup of $K^*/k^*\otimes \Z_{(\ell)}$.

\

{\em Step 10.} By Proposition~\ref{prop:KKK},
$\sK^*$ is a multiplicative group of a field so that both 
$K$ and $L$ are finite purely inseparable extensions of this 
field.
It remains to insure that the additive structure on $\sK^*$ 
is intrinsically defined.

\

{\em Step 11.} By Theorem~\ref{thm:skk} and  
Proposition~\ref{prop:unique-p}, the field is uniquely determined
by the partial projective structure.

\

{\em Step 12.} Lemma~\ref{lemm:unnn}  and Corollary~\ref{coro:pro-uni}
give a Galois-theoretic characterization of 
generating elements and primary lines in $\sK^*/k^*$.
Proposition~\ref{prop:gener} and 
Example~\ref{exam:basic2} show that these 
define a (unique) partial projective structure on $\sK^*/k^*$
(in particular, the projective structures induced by 
$\P(K)$ and $\P(L)$
coincide).  

\

{\em Step 13.} If follows $K/\sK$ and $L/\sK$ are finite purely inseparable
extensions of the {\em same} field. This concludes the proof 
of Theorem~\ref{thm:main}.

\bibliographystyle{smfplain}
\bibliography{recon}

\end{document}